\documentclass[10pt]{amsart}
\usepackage{amscd}
\usepackage{amssymb}
\usepackage[all]{xy}
\newcommand{\choosegraphics}[1]{#1}

\newtheorem{thm}[equation]{Theorem}
\newtheorem{cor}[equation]{Corollary}
\newtheorem{prop}[equation]{Proposition}
\newtheorem{lem}[equation]{Lemma}
\theoremstyle{definition}
\newtheorem{dfn}[equation]{Definition}
\newtheorem{rem}[equation]{Remark}
\newtheorem{exa}[equation]{Example}

\newtheorem{que}[equation]{Question}
\numberwithin{equation}{section}

\newcommand{\iso}{\stackrel{\simeq}{\rightarrow}}
\newcommand{\inj}{\hookrightarrow}
\newcommand{\surj}{\twoheadrightarrow}
\newcommand{\opn}{\operatorname}
\newcommand{\cat}[1]{\operatorname{\mathsf{#1}}}

\newcommand{\rmitem}[1]{\item[\text{\textup{(#1)}}]}

\newcommand{\mfrak}[1]{\mathfrak{#1}}
\newcommand{\mcal}[1]{\mathcal{#1}}
\newcommand{\msf}[1]{\mathsf{#1}}
\newcommand{\mbf}[1]{\mathbf{#1}}
\newcommand{\mrm}[1]{\mathrm{#1}}
\newcommand{\mbb}[1]{\mathbb{#1}}

\newcommand{\tup}[1]{\textup{#1}}
\newcommand{\bsym}[1]{\boldsymbol{#1}}
\newcommand{\bra}[1]{\langle #1 \rangle}

\begin{document}
\title{Residue Complexes over Noncommutative Rings}
\author{Amnon Yekutieli and James J.\ Zhang }
\address{Department of Mathematics,
Ben Gurion University, Be'er Sheva 84105, ISRAEL}
\email{amyekut@math.bgu.ac.il}
\urladdr{http://www.math.bgu.ac.il/$\sim$amyekut}
\address{Department of Mathematics, University of
Washington, Seattle, WA 98195, USA}
\email{zhang@math.washington.edu}
\date{1 June 2002}
\thanks{Both authors were supported by the US-Israel Binational
Science Foundation. The second author was also supported by the 
NSF, a Sloan Research Fellowship and the Royalty Research Fund of 
the University of Washington}

\begin{abstract}
Residue complexes were introduced by Grothendieck in algebraic 
geometry. These are canonical complexes of injective modules that 
enjoy remarkable functorial properties (traces). 

In this paper we study residue complexes over noncommutative 
rings. These objects have a more intricate structure than in the 
commutative case, since they are complexes of bimodules. We 
develop methods to prove uniqueness, existence and functoriality 
of residue complexes.

For a polynomial identity algebra over a field (admitting a 
noetherian connected filtration) we prove existence of the residue 
complex and describe its structure in detail.
\end{abstract}

\maketitle


\setcounter{section}{-1}
\section{introduction}

\subsection{Motivation: A Realization of the Geometry of a 
Noncommutative Ring} 
For a commutative ring $A$ it is clear (since Grothendieck) what 
is the geometric object associated to $A$: the locally 
ringed space $\opn{Spec} A$. 
However if $A$ is noncommutative this question becomes pretty 
elusive. One possibility is to consider the set $\opn{Spec} A$ of 
two-sided prime (or maybe primitive) ideals of $A$. Another 
possibility is to choose a side -- say left -- and to consider the 
category $\cat{Mod} A$ of left $A$-modules (or some related 
construction) as a kind of geometric object. Both these options 
are used very effectively in various contexts; but neither is 
completely satisfactory. A common (genuine) obstacle is the 
difficulty of localizing noncommutative rings. A ``classical'' 
account of the subject can be found in \cite{MR}; recent 
developments are described in \cite{SV} and its references.

In this paper we try another point of view. Taking our cue from 
commutative algebraic geometry, we try to construct a global 
algebraic object -- the {\em residue complex $\mcal{K}_{A}$} -- which 
encodes much of the geometric information of $A$. 

Let us first examine an easy case  
which can explain where we are heading. Suppose $\mbb{K}$ is a 
field and $A$ is a finite $\mbb{K}$-algebra. If $A$ is commutative 
then $\opn{Spec} A$ is a finite set. The injective module 
$A^{*} := \opn{Hom}_{\mbb{K}}(A, \mbb{K})$
is a direct sum of indecomposable modules, each summand 
corresponding a point in $\opn{Spec} A$.

On the other hand if $A$ is noncommutative the geometric object 
associated to it is a finite quiver $\vec{\Delta}$. The vertex set 
of $\vec{\Delta}$ is $\opn{Spec} A$, and the arrows (links) are 
determined by the bimodule decomposition of 
$\mfrak{r} / \mfrak{r}^{2}$, where $\mfrak{r}$ is the 
Jacobson radical. The connected components of $\vec{\Delta}$ are 
called cliques. Here is the corresponding module-theoretic 
interpretation: the vertices of $\vec{\Delta}$ are the isomorphism 
classes of indecomposable summands of $A^{*}$ as left module, 
the cliques are the indecomposable summands of $A^{*}$ as 
bimodule, and the arrows in $\vec{\Delta}$ represent irreducible 
homomorphisms between vertices. Finally if $A \to B$ is a finite 
homomorphism then there is an $A$-bimodule homomorphism 
$\opn{Tr}_{B / A} : B^{*} \to A^{*}$.

The point of view we adopt in this paper is that for some
infinite noncommutative $\mbb{K}$-algebras $A$ the module-theoretic 
interpretation of the geometry of $A$, as stated above, should also 
make sense. The generalization of the bimodule 
$A^{*}$ is the residue complex $\mcal{K}_{A}$. The additional data 
(not occurring in finite algebras) is that of specialization, 
which should be carried by the coboundary operator of $\mcal{K}_{A}$.

There are certain cases in which we know this plan works. 
For commutative rings this is Grothendieck's theory of residual 
complexes, worked out in \cite{RD}, and reviewed in Subsection 0.2 
below. If $A$ is finite over its center $\mrm{Z}(A)$ then the cliques 
of $A$ biject to $\opn{Spec} \mrm{Z}(A)$, and hence the geometry 
of $A$ is understood; and the residue complex is
$\mcal{K}_{A} = \opn{Hom}_{\mrm{Z}(A)}(A, \mcal{K}_{\mrm{Z}(A)})$.

If $A$ is a twisted homogeneous coordinate ring 
of a projective variety $X$ (with automorphism $\sigma$ and 
$\sigma$-ample line bundle $\mcal{L}$) we know the graded residue 
complex $\mcal{K}_{A}$ exists (see \cite{Ye1}). 
Here the indecomposable graded left module summands 
of $\mcal{K}^{-q - 1}_{A}$, $0 \leq q \leq \opn{dim} X$, 
are indexed by the points of $X$ of 
dimension $q$; and the indecomposable graded bimodule summands are the
$\sigma$-orbits of these points. 
A similar phenomenon (for $q = 0,1)$ occurs when $A$ is a 
$3$-dimensional Sklyanin algebra (see \cite{Ye2}). 

In Subsection 0.3 we give a brief explanation of the 
noncommutative residue complex, and state the main results of our 
paper.

\subsection{R\'esum\'e: Residue Complexes in Algebraic Geometry}
Residue complexes in (commutative) algebraic geometry were introduced
by Grothendieck \cite{RD}. 
Suppose $\mbb{K}$ is a field and $X$ is a finite type $\mbb{K}$-scheme.
The residue complex of $X$ is a bounded complex 
$\mcal{K}_{X}$ of quasi-coherent sheaves with some remarkable 
properties. First each of the $\mcal{O}_{X}$-modules 
$\mcal{K}^{-q}_{X}$ is injective, and the functor
$\mcal{H}om_{\mcal{O}_{X}}(-, \mcal{K}_{X})$
is a duality of the bounded derived category with coherent 
cohomology
$\msf{D}^{\mrm{b}}_{\mrm{c}}(\cat{Mod} \mcal{O}_{X})$.
Next, if $f: X \to Y$ is a proper morphism then there is a 
nondegenerate trace map (an actual homomorphism of complexes)
$\opn{Tr}_{f}: f_{*} \mcal{K}_{X} \to \mcal{K}_{Y}$. 
Finally if $X$ is smooth of dimension $n$ over $\mbb{K}$ then there 
is a canonical quasi-isomorphism
$\Omega^{n}_{X / \mbb{K}}[n] \to \mcal{K}_{X}$. 

In \cite{RD} the residue complex $\mcal{K}_{X}$ is closely related 
to the twisted inverse image functor. Indeed, if we denote by 
$\pi_{X}: X \to \opn{Spec} \mbb{K}$ the structural morphism, then 
the twisted inverse image 
$\pi_{X}^{!} \mbb{K} \in 
\msf{D}^{\mrm{b}}_{\mrm{c}}(\cat{Mod} \mcal{O}_{X})$ 
is a dualizing complex. There is a trace 
$\opn{Tr}_{f}: \mrm{R} f_{*} \pi_{X}^{!} \mbb{K} \to
\pi_{Y}^{!} \mbb{K}$
for a proper morphism $f: X \to Y$, and an isomorphism
$\Omega^{n}_{X / \mbb{K}}[n] \to \pi_{X}^{!} \mbb{K}$
for $X$ smooth. 

The filtration of $\cat{Mod} \mcal{O}_{X}$ by dimension of support 
(niveau filtration) gives rise to the {\em Cousin functor} $\mrm{E}$.
For a complex $\mcal{M}$ the Cousin complex $\mrm{E} \mcal{M}$ 
is the row $q = 0$ in the $E_{1}$ page of the niveau spectral 
sequence
$E_{1}^{p, q} \Rightarrow \mrm{H}^{p + q} \mcal{M}$.
In this way one obtains a functor
$\mrm{E}: \msf{D}^{+}(\cat{Mod} \mcal{O}_{X}) \to 
\msf{C}^{+}(\cat{Mod} \mcal{O}_{X})$
where the latter is the (abelian) category of complexes. 
By definition the residue complex is 
$\mcal{K}_{X} := \mrm{E} \pi_{X}^{!} \mbb{K}$,
and there is a canonical isomorphism 
$\pi_{X}^{!} \mbb{K} \cong \mcal{K}_{X}$
in the derived category. Explicit constructions of the residue 
complex also exist; cf.\ \cite{Ye3} and its references.

Here is what this means for affine schemes. If we consider a 
commutative finitely generated $\mbb{K}$-algebra $A$, and 
$X := \opn{Spec} A$, then 
$\mcal{K}_{A} := \Gamma(X, \mcal{K}_{X})$
is a bounded complex of injective $A$-modules. 
For any integer $q$ there is a decomposition
$\mcal{K}_{A}^{-q} \cong \bigoplus J(\mfrak{p})$,
where $\mfrak{p}$ runs over the prime ideals such that 
$\opn{dim} A / \mfrak{p} = q$, and $J(\mfrak{p})$ is the 
injective hull of $A / \mfrak{p}$. 
The map 
$J(\mfrak{p}) \inj \mcal{K}_{A}^{-q} \to \mcal{K}_{A}^{-q + 1}
\surj J(\mfrak{q})$
is nonzero precisely when $\mfrak{p} \subset \mfrak{q}$.
Moreover, $\mcal{K}_{A}$ is dualizing, in the sense that the 
functor 
$\opn{Hom}_{A}(-, \mcal{K}_{A})$
is a duality of the bounded derived category with finite 
cohomologies
$\msf{D}^{\mrm{b}}_{\mrm{f}}(\cat{Mod} A)$.
If $A \to B$ is a finite homomorphism then there is a nondegenerate 
trace map
$\opn{Tr}_{B / A}: \mcal{K}_{B} \to \mcal{K}_{A}$. 
And if $A$ is smooth of relative dimension $n$ then
\[ 0 \to \Omega^{n}_{A / \mbb{K}} \to \mcal{K}^{-n}_{A} \to \cdots
\mcal{K}^{0}_{A} \to 0 \]
is a minimal injective resolution.

\subsection{Statement of Main Results}
In the present paper we study a noncommutative version of the 
above. Now $A$ is an associative, unital, noetherian, affine 
(i.e.\ finitely generated) $\mbb{K}$-algebra, not necessarily 
commutative. We denote by $\cat{Mod} A$ the category of left 
$A$-modules and by $A^{\mrm{op}}$ the opposite algebra.

A dualizing complex over the algebra $A$ is, roughly speaking,
a complex $R$ of bimodules,  
such that the two derived functors $\opn{RHom}_{A}(-, R)$ and
$\opn{RHom}_{A^{\mrm{op}}}(-, R)$ induce a duality between 
$\msf{D}^{\mrm{b}}_{\mrm{f}}(\cat{Mod} A)$ and 
$\msf{D}^{\mrm{b}}_{\mrm{f}}(\cat{Mod} A^{\mrm{op}})$.
The full definition of this, as well as of other important notions,
are included in the body of the paper. Dualizing 
complexes over noncommutative rings have various applications, for 
instance in ring theory (see \cite{YZ2}), 
representation theory (see \cite{Ye5}, \cite{EG} and \cite{BGK}),  
and even theoretical physics (see \cite{KKO}).

The twisted inverse image $\pi^{!}_{X} \mbb{K}$ of the commutative 
picture is generalized to the {\em rigid dualizing complex} $R$, 
as defined by Van den Bergh \cite{VdB1}. Indeed if $A$ is 
commutative and $X = \opn{Spec} A$ then 
$R := \mrm{R} \Gamma(X, \pi_{X}^{!} \mbb{K})$ is a rigid dualizing 
complex. For noncommutative $A$ we know that a rigid dualizing 
complex $R_{A}$ (if exists) is unique, and for a finite homomorphism 
$A \to B$ there is at most one rigid trace
$\opn{Tr}_{B / A}: R_{B} \to R_{A}$. 

It is known \cite{YZ2} that if $R$ is an {\em Auslander dualizing 
complex} then the canonical dimension associated to $R$, namely
$\opn{Cdim} M := - \opn{inf} \{ q \mid \opn{Ext}^{q}_{A}(M, R) 
\neq 0 \}$
for a finite module $M$, is an exact dimension function. 

The residue complex of $A$ is by definition a rigid Auslander 
dualizing complex $\mcal{K}_{A}$, consisting of bimodules 
$\mcal{K}_{A}^{-q}$ which are injective, and pure of dimension 
$q$ with respect to $\opn{Cdim}$, on both sides. Again, if $A$ is 
commutative then this definition is equivalent to that of 
\cite{RD}. 

The Cousin functor is available in the noncommutative situation 
too. Assume we have a rigid Auslander dualizing complex $R_{A}$. 
The canonical dimension $\opn{Cdim}$ gives a filtration of 
$\cat{Mod} A$ by ``dimension of support'', and just like in the 
commutative case we obtain a Cousin functor 
$\mrm{E}: \msf{D}^{+}(\cat{Mod} A \otimes A^{\mrm{op}}) \to 
\msf{C}^{+}(\cat{Mod} A \otimes A^{\mrm{op}})$.
However usually $\mrm{E} R_{A}$ will not be a residue complex!

The first main result gives a sufficient condition for the existence 
of a residue complex (it is not hard to see that this condition
is also necessary). We say $R_{A}$ has a pure
minimal injective resolution on the left if in the minimal 
injective resolution $R_{A} \to I$ in $\msf{C}^{+}(\cat{Mod} A)$ 
each $I^{-q}$ is pure of $\opn{Cdim} = q$; likewise on the right.

\begin{thm} \label{thm0.1}
Suppose $A$ is a noetherian $\mbb{K}$-algebra and $R_{A}$ is an 
Auslander rigid dualizing complex over $A$. Assume $R_{A}$ has pure 
minimal injective resolutions on both sides. Then 
$\mcal{K}_{A} := \mrm{E} R_{A}$
is a residue complex.
\end{thm}

This result included in Theorem \ref{thm4.1} in the body of the 
paper. We also have a result guaranteeing the existence of a 
trace between residue complexes (it is part of Theorem \ref{thm5.1}). 
One calls a ring homomorphism $A \to B$ a finite 
centralizing homomorphism if 
$B = \sum A b_{i}$ where the $b_{i}$ are finitely many elements 
of $B$ that commute with every $a \in A$. 

\begin{thm} \label{thm0.4}
Let $A \to B$ be a finite centralizing homomorphism between
noetherian $\mbb{K}$-algebras. Suppose the two conditions below
hold.
\begin{enumerate}
\rmitem{i} There are  rigid dualizing complexes
$R_{A}$ and $R_{B}$ and the rigid trace morphism
$\opn{Tr}_{B / A} : R_{B} \to R_{A}$ exists.
\rmitem{ii} $R_{A}$ is an Auslander dualizing complex and it has
pure minimal injective resolutions on both sides.
\end{enumerate}
Let $\mcal{K}_{A} := \mrm{E} R_{A}$ be the residue complex of $A$ 
\tup{(}cf.\ Theorem \tup{\ref{thm0.1})}.
Then $\mcal{K}_{B} := \mrm{E} R_{B}$ is the residue complex of 
$B$. The homomorphism of complexes
$\mrm{E}(\opn{Tr}_{B / A}): \mcal{K}_{B} \to \mcal{K}_{A}$
is a rigid trace, and it induces an isomorphism of complexes of 
$A$-bimodules
\[ \mcal{K}_{B} \cong \opn{Hom}_{A}(B, \mcal{K}_{A}) \cong
\opn{Hom}_{A^{\mrm{op}}}(B, \mcal{K}_{A})  . \]
\end{thm}

In Subsection 0.1 we listed a few classes of algebras for which 
residue complexes were previously known to exist. More examples appear 
in Section 5 of the paper (e.g.\ the first Weyl algebra, the universal
enveloping algebra of a nilpotent $3$-dimensional Lie algebra). In 
the remainder of this subsection we discuss our results for the class 
of polynomial identity (PI) algebras. 

We remind that a PI ring $A$ is one that 
satisfies some polynomial identity
$f(x_{1}, \ldots, x_{n}) = 0$, and hence is close to being 
commutative (a commutative ring satisfies the identity 
$x_1 x_2 - x_2 x_1 = 0$). We have a quite detailed knowledge of the 
residue complex of a PI algebra A, assuming it admits some
noetherian connected filtration. 
A a noetherian connected filtration on the algebra $A$ 
is a filtration $\{ F_{n} A \}$ such that the graded algebra 
$\opn{gr}^{F} A$ is a noetherian connected graded 
$\mbb{K}$-algebra. Most known examples of noetherian affine PI
algebras admit noetherian connected filtrations, but there are 
counterexamples (see \cite{SZ}). 

\begin{thm} \label{thm0.2}
Let $A$ be an affine noetherian PI algebra admitting a noetherian 
connected filtration.
\begin{enumerate}
\item
$A$ has a residue complex $\mcal{K}_{A}$.
\item
Let  $B = A / \mfrak{a}$ be a quotient algebra. Then $B$ has a residue
complex $\mcal{K}_{B}$, there is a rigid trace
$\opn{Tr}_{B / A} : \mcal{K}_{B} \to \mcal{K}_{A}$
that is an actual homomorphism of complexes of bimodules, and 
$\opn{Tr}_{B / A} $ induces an isomorphism
\[ \mcal{K}_{B} \cong \opn{Hom}_{A}(B, \mcal{K}_{A}) =
\opn{Hom}_{A^{\mrm{op}}}(B, \mcal{K}_{A}) \subset \mcal{K}_{A} . \]
\end{enumerate}
\end{thm}

This is Theorem \ref{thm6.1} in the body of the paper. 

The next theorem describes the structure of the residue complex 
$\mcal{K}_{A}$ of a PI algebra. 

Recall that the prime spectrum $\opn{Spec} A$ is a disjoint union 
of cliques. For any clique $Z$ we denote by $A_{S(Z)}$ the 
localization at $Z$, and for a module $M$ we let $\Gamma_{Z} M$ 
be the submodule supported on $Z$. 
We say a clique $Z_{1}$ is a specialization of a clique $Z_{0}$ if 
there are prime ideals $\mfrak{p}_{i} \in Z_{i}$ with
$\mfrak{p}_{0} \subset \mfrak{p}_{1}$. 

The $q$-skeleton of $\opn{Spec} A$ is the set of 
prime ideals $\mfrak{p}$ such that $\opn{Cdim} A / \mfrak{p} = q$.
It is a union of cliques. 

For any prime ideal $\mfrak{p}$ we let $J(\mfrak{p})$ be the 
indecomposable injective $A$-module with associated prime 
$\mfrak{p}$, and $\mrm{r}(\mfrak{p})$ is the Goldie rank of 
$A / \mfrak{p}$.

\begin{thm} \label{thm0.3}
Let $A$ be a PI $\mbb{K}$-algebra admitting a 
noetherian connected filtration, and let $\mcal{K}_{A}$ be its 
residue complex.
\begin{enumerate}
\item For every $q$ there is a canonical $A$-bimodule 
decomposition
\[ \mcal{K}^{-q}_{A} = \bigoplus_{Z} \Gamma_{Z} \mcal{K}^{-q}_{A} \]
where $Z$ runs over the cliques in the $q$-skeleton of
$\opn{Spec} A$.
\item Fix one clique $Z$ in the $q$-skeleton of $\opn{Spec} A$.
Then $\Gamma_{Z} \mcal{K}^{-q}_{A}$ is an indecomposable
$A$-bimodule.
\item $\Gamma_{Z} \mcal{K}^{-q}_{A}$
is an injective left $A_{S(Z)}$-module, 
and its socle is the essential submodule
\[ \bigoplus_{\mfrak{p} \in Z}  \mcal{K}^{-q}_{A / \mfrak{p}} \cong
\bigoplus_{\mfrak{p} \in Z}
\opn{Hom}_{A}(A / \mfrak{p}, \mcal{K}^{-q}_{A})
\subset \Gamma_{Z} \mcal{K}_{A}^{-q} . \]
\item There is a \tup{(}noncanonical\tup{)} decomposition of left 
$A_{S(Z)}$-modules
\[ \Gamma_{Z} \mcal{K}^{-q}_{A} 
\cong \bigoplus_{\mfrak{p} \in Z} 
J_{A}(\mfrak{p})^{\mrm{r}(\mfrak{p})} . \]
\item Suppose $Z_{i}$ is a clique in the $(q - i)$-skeleton of
$\opn{Spec} A$, for $i = 0, 1$.
Then $Z_{1}$ is a specialization of $Z_{0}$ iff the
composed homomorphism
\[ \Gamma_{Z_{0}} \mcal{K}^{-q}_{A} \inj
\mcal{K}^{-q}_{A}
\to \mcal{K}^{-q + 1}_{A} \surj \Gamma_{Z_{1}}
\mcal{K}^{-q + 1}_{A} \]
is nonzero.
\end{enumerate}
\end{thm}

This theorem is repeated as Theorem \ref{thm6.2} in the body of 
the paper. 

Observe that part (4) of the theorem says that the prime 
spectrum $\opn{Spec} A$ is encoded in the left module decomposition 
of the complex $\mcal{K}_{A}$. By the left-right symmetry 
(replacing $A$ with $A^{\mrm{op}}$) the same is true for the right 
module decomposition of $\mcal{K}_{A}$. 
Parts (1) and (2) imply that the cliques in $\opn{Spec} A$ are 
encoded in the bimodule decomposition of $\mcal{K}_{A}$. Part (5) 
says that specializations are encoded in the coboundary operator 
of $\mcal{K}_{A}$.

We end this Subsection with a disclaimer.
If the algebra $A$ is ``too noncommutative'' then it will not 
have a residue complex (for instance $A = \mrm{U}(\mfrak{sl}_{2})$). 
Thus the scope of the theory of residue complexes is necessarily 
limited. Our upcoming paper \cite{YZ3} presents an alternative 
approach to address precisely this issue.

\subsection{Outline of the Paper} \mbox{ }

\medskip \noindent
\textbf{Section 1.}
We begin by recalling the notions of localizing subcategories and 
torsion functors in the module category
$\cat{Mod} A$ of a ring $A$. Given a localizing subcategory 
$\msf{M} \subset \cat{Mod} A$ we consider the derived functor 
$\mrm{R} \Gamma_{\msf{M}}$, and the cohomology with support in 
$\msf{M}$, namely
$\mrm{H}^{q}_{\msf{M}} := \mrm{H}^{q} \mrm{R} \Gamma_{\msf{M}}$.
We recall what is an exact dimension function, and relate it to 
localizing subcategories. $\msf{M}$-flasque modules are defined.
The main result here is Theorem \ref{thm1.1}, dealing with 
cohomology with supports for bimodules.
 
\medskip \noindent
\textbf{Section 2.}
Given a filtration $\msf{M} = \{ \msf{M}^{p} \}$ 
of $\cat{Mod} A$ by localizing subcategories we can define the 
Cousin functor 
$\mrm{E}_{\msf{M}}: \msf{D}^{+}(\cat{Mod} A) \to 
\msf{C}^{+}(\cat{Mod} A)$. 
The main result here is Theorem \ref{thm2.1} which provides a 
sufficient condition for a complex $M$ to be isomorphic to its 
Cousin complex
$\mrm{E}_{\msf{M}} M$ in $\msf{D}^{+}(\cat{Mod} A)$.

\medskip \noindent
\textbf{Section 3.}
The definitions of rigid dualizing complex and rigid trace are
recalled. We show that the rigid dualizing complex is compatible 
with central localization. 

\medskip \noindent
\textbf{Section 4.}
Here we look at residual complexes, which are Auslander 
dualizing complexes consisting of bimodules that are injective 
and pure on both sides. We prove Theorem \ref{thm4.1}, which is 
the essential ingredient of Theorem \ref{thm0.1} above. Also we 
prove Theorem \ref{thm4.2}, asserting the existence of a residual 
complex over an FBN algebra $A$ with an Auslander dualizing 
complex satisfying a certain symmetry condition.

\medskip \noindent
\textbf{Section 5.}
A residue complex is a residual complex that is also rigid. The 
main result in this section is Theorem \ref{thm5.1} which relates 
the rigid trace to residue complexes. As a corollary we deduce 
that a residue complex $\mcal{K}_{A}$ over an algebra $A$ is 
unique up to a unique isomorphism of complexes (Corollary 
\ref{cor5.0}). We present examples of algebras with residue 
complexes. Also we explain what our results 
mean for commutative algebras.

\medskip \noindent
\textbf{Section 6.}
Besides proving Theorems \ref{thm0.2} and \ref{thm0.3}, we also 
prove that for a prime PI algebra $A$ of dimension $n$
the generic component $\mcal{K}_{A}^{-n}$ is untwisted; in fact it 
is isomorphic as bimodule to the ring of fractions $Q$. Several 
examples are studied.

\bigskip \noindent
\textbf{Acknowledgements.}
We wish to thank our teacher M. Artin who brought the question of 
residue complexes to our attention over a decade ago. Many thanks to 
M. Van den Bergh who shared with us many good ideas, and pointed 
out a mistake in an earlier version of the paper. Also thanks to 
K. Goodearl, C. Ingalls, J. Lipman, S.P. Smith, J.T. Stafford and 
Q.S. Wu for suggestions and discussions.

\section{Cohomology with Support in a Localizing Subcategory}

In algebraic geometry, given a scheme $X$ and a closed subset
$Z \subset X$, the functor $\Gamma_{Z}$ is defined: for any sheaf
$\mcal{M}$, $\Gamma_{Z} \mcal{M} \subset \mcal{M}$ is the subsheaf
of sections supported on $Z$. The derived functors
$\mrm{H}^{q} \mrm{R} \Gamma_{Z} \mcal{M} = \mcal{H}^{q}_{Z} \mcal{M}$
are called the sheaves of cohomologies of $\mcal{M}$ with support
in $Z$. More generally one can take a family of supports
$\mcal{Z}$, which is a family of closed sets satisfying suitable
conditions (e.g.\
$\mcal{Z}_{q} = \{ Z \text{ closed, } \opn{dim} Z \leq q \}$).

In this section we consider an analogous construction replacing
the scheme $X$ with the category $\cat{Mod} A$ of left modules
over a ring $A$.
The role of family of supports is played by a localizing subcategory
$\cat{M} \subset \cat{Mod} A$. This idea already appeared in
\cite{Ye2}, but here we expand the method significantly.
With minor modifications the contents of this section and the
next one will apply to any noetherian quasi-scheme $X$ (in the
sense of \cite{VdB2}).

We begin with a quick review of Gabriel's theory of torsion,
following \cite[Chapter VI]{Ste}, but using notation suitable for
our purposes. Fix a ring $A$.
A {\em left exact radical} (or {\it torsion functor}) is an additive
functor $\Gamma : \cat{Mod} A \to \cat{Mod} A$,
which is a subfunctor of the identity functor
$\bsym{1}_{\cat{Mod} A}$,
left exact, and $\Gamma(M / \Gamma M) = 0$ for any
$M \in \cat{Mod} A$. It follows that
$\Gamma \Gamma M  = \Gamma M$, and if $N \subset M$ then
$\Gamma N = N \cap \Gamma M$.

A {\em hereditary torsion class} is a class of objects
$\cat{M} \subset \cat{Mod} A$ closed under subobjects, quotients,
extensions and infinite direct sums. The full subcategory
$\cat{M}$ is a {\em localizing subcategory}.
Given a left exact radical $\Gamma$, the subcategory
\[ \cat{M}_{\Gamma} := \{ M \mid \Gamma M = M \} \]
is localizing. Conversely, given a localizing subcategory
$\cat{M}$, the functor
\[ M \mapsto \Gamma_{\cat{M}} M :=
\{ m \in M \mid A \cdot m \in \cat{M} \} \]
is a left exact radical. One has
$\Gamma_{\msf{M}_{\Gamma}} = \Gamma$ and 
$\msf{M}_{\Gamma_{\msf{M}}} = \msf{M}$.

A third equivalent notion is that of {\em left Gabriel topology}
(or {\em filter}) in $A$, which is a set $\mfrak{F}$ of left ideals
of $A$ satisfying some axioms (that we shall not list here).
Given a localizing subcategory
$\cat{M} \subset \cat{Mod} A$,
the set of left ideals
\[ \mfrak{F}_{\cat{M}} := \{ \mfrak{a} \subset A \text{ left ideal }
\mid A / \mfrak{a} \in \cat{M} \} \]
is a left Gabriel topology, and any left Gabriel topology arises this
way. On the other hand given a left Gabriel
topology $\mfrak{F}$ the functor
\begin{equation} \label{eqn1.1}
\Gamma_{\mfrak{F}} := \underset{\mfrak{a} \in \mfrak{F}}{\varinjlim}
\opn{Hom}_{A}(A / \mfrak{a}, - ) ,
\end{equation}
where $\mfrak{F}$ is partially ordered inclusion,
is a left exact radical.

Below are some examples of localizing subcategories.

\begin{exa} \label{exa1.1}
Let $Z$ be a set of two-sided ideals of a ring $A$, each of which 
is finitely generated as left ideal. Then the set
\begin{equation} \label{eqn1.2}
\mfrak{F}_{Z} := \{ \text{ left ideals }\mfrak{a} \subset A  \mid
\mfrak{m}_{1} \cdots \mfrak{m}_{n} \subset \mfrak{a}
\text{ for some } \mfrak{m}_{1}, \ldots, \mfrak{m}_{n} \in Z \} 
\end{equation}
is a left Gabriel topology (cf.\ \cite[Proposition VI.6.10]{Ste}).
The corresponding torsion functor is denoted $\Gamma_{Z}$ and the 
localizing subcategory is $\msf{M}_{Z}$. 
If $Z = \{ \mfrak{m} \}$ we also write
$\Gamma_{\mfrak{m}}$ and $\msf{M}_{\mfrak{m}}$.
If $A$ is commutative then $\mfrak{F}_{\mfrak{m}}$ is the usual
$\mfrak{m}$-adic topology, and $\Gamma_{\mfrak{m}} M$
is the submodule of elements supported on $\opn{Spec} A / \mfrak{m}$.
\end{exa}


Keeping this example in mind, in the general situation of a localizing
subcategory $\cat{M}$ we call $\Gamma_{\cat{M}} M$ the submodule of
elements supported on $\cat{M}$.

A localizing subcategory $\cat{M}$ is called {\em stable} if
whenever $M \in \cat{M}$  and $M \subset N$ is an essential
submodule then also $N \in \cat{M}$.

\begin{exa} \label{exa1.4}
Suppose $A$ is left noetherian. If the ideal $\mfrak{a}$
has the left Artin-Rees property (e.g.\ when $\mfrak{a}$ is generated
by normalizing elements) then the localizing subcategory
$\msf{M}_{\mfrak{a}}$ is stable. See \cite[Theorem 2.2 and Proposition
2.6]{MR}. More generally, if the set $Z$ of ideals has
the Artin-Rees property then $\msf{M}_{Z}$ is stable, see
\cite[Proposition 2.9]{BM}.
\end{exa}

In this paper the most important examples of localizing
subcategories arise from dimension functions.

\begin{dfn}
Let $\cat{M}$ be an abelian category. 
An {\em exact dimension function} on $\cat{M}$ is a function
$\opn{dim}: \cat{M} \to \{ -\infty \} \cup \mbb{R}
\cup \{ \text{infinite ordinals}\}$,
satisfying the following axioms:
\begin{enumerate}
\rmitem{i} $\opn{dim} 0 = -\infty$.
\rmitem{ii} For every short exact sequence
$0 \to M' \to M \to M'' \to 0$ one has
$\opn{dim} M = \max \{ \opn{dim} M', \opn{dim} M'' \}$.
\rmitem{iii} If $M = \lim_{\alpha \to} M_{\alpha}$
and each $M_{\alpha} \to M$ is an injection then
$\opn{dim} M  = \opn{sup} \{ \opn{dim} M_{\alpha} \}$.
\end{enumerate}
\end{dfn}

When $\cat{M} = \cat{Mod} A$ for a left noetherian ring $A$, 
often a dimension function will satisfy a further axiom.
Recall that if $\mfrak{p}$ is a prime ideal then 
an $A / \mfrak{p}$-module $M$ is called {\em torsion} if for any 
$m \in M$ there is a regular element $a \in A / \mfrak{p}$ 
such that $a m = 0$. 

\begin{dfn} \label{dfn1.1}
Let $A$ be a left noetherian ring. A {\em spectral exact dimension 
function} on $\cat{Mod} A$ is an exact dimension function
$\opn{dim}$ satisfying the extra axiom
\begin{enumerate}
\rmitem{iv} If $\mfrak{p} M = 0$ for some prime ideal $\mfrak{p}$,
and $M$ is a torsion $A / \mfrak{p}$-module, then
$\opn{dim} M \leq \opn{dim} A / \mfrak{p} - 1$.
\end{enumerate}
\end{dfn}

In this paper the dimension functions will all take values in 
$\{ -\infty \} \cup \mbb{Z}$.

\begin{rem}
The definition of spectral exact dimension function is 
standard in ring theory, although usually one restricts to the 
subcategory $\cat{Mod}_{\mrm{f}} A$ of finite (i.e.\ finitely 
generated) modules, where condition (iii) becomes trivial.
Cf.\ \cite[Section 6.8.4]{MR}. 
\end{rem}

\begin{exa} \label{exa1.2}
Let $\opn{dim}$ be an exact dimension function on $\cat{Mod} A$.
For an integer $q$ let
\[ \msf{M}_{q}(\opn{dim})
:= \{ M \in \cat{Mod} A \mid \opn{dim} M \leq q \} . \]
Then $\msf{M}_{q}(\opn{dim})$ is a localizing subcategory.
\end{exa}

Here is a different kind of localizing subcategory.

\begin{exa} \label{exa1.3}
Given a left denominator set $S \subset A$ 
(cf.\ \cite[Paragraph 2.1.13]{MR}) we define
\[ \mfrak{F}_{S} := \{ \mfrak{a} \subset A \text{ left ideal } \mid
\mfrak{a} \cap S \neq \emptyset \} . \]
According to \cite[Section II.3 and Example in Section VI.9]{Ste} 
this is a left Gabriel topology. We denote the
localizing subcategory by $\msf{M}_{S}$. Letting $A_{S} = S^{-1}
A$ be the left ring of fractions, for every module $M$ one has an
exact sequence
\[ 0 \to \Gamma_{\msf{M}_{S}} M \to M \to A_{S} \otimes_{A} M . \]
\end{exa}

Now we want to pass to derived categories.
Let $\msf{D}(\cat{Mod} A)$ be the derived category of $A$-modules,
and let $\msf{D}^{+}(\cat{Mod} A)$ be the full subcategory
of bounded below complexes. As usual $\msf{C}(\cat{Mod} A)$ denotes
the abelian category of complexes. Our references are
\cite[Section I]{RD} and \cite[Section I]{KS}.

Suppose $M \in \msf{C}^{+}(\cat{Mod} A)$. By an {\em injective
resolution of $M$ in} $\msf{C}^{+}(\cat{Mod} A)$ we mean a
quasi-isomorphism $M \to I$ in $\msf{C}^{+}(\cat{Mod} A)$ with
each $I^{q}$ an injective module.

\begin{lem}
Let $\cat{M}$ be a localizing subcategory of $\cat{Mod} A$.
Then there is a right derived functor
\[ \mrm{R} \Gamma_{\cat{M}} : \cat{D}^{+}(\cat{Mod} A) \to
\cat{D}^{+}(\cat{Mod} A) . \]
\end{lem}

\begin{proof}
Given $M \in \cat{D}^{+}(\cat{Mod} A)$ take any injective
resolution $M \to I$ in $\cat{C}^{+}(\cat{Mod} A)$, and let
$\mrm{R} \Gamma_{\cat{M}} M := \Gamma_{\cat{M}} I$
(cf.\ \cite[Theorem I.5.1]{RD}).
\end{proof}

Note that
$\mrm{R} \Gamma_{\cat{M}} M \in \cat{D}^{+}_{\cat{M}}(\cat{Mod} A)$,
the full triangulated subcategory whose objects are
complexes with cohomology in $\cat{M}$.

\begin{rem}
One can define $\mrm{R} \Gamma_{\cat{M}} M$ for an unbounded complex
as $\mrm{R} \Gamma_{\cat{M}} M := \Gamma_{\cat{M}} I$ where
$M \to I$ is a quasi-isomorphism to a K-injective complex $I$,
cf.\ \cite{Sp}.
\end{rem}

\begin{dfn}
The $q$th cohomology of $M$ with support in $\msf{M}$ is defined
to be
$\mrm{H}^{q}_{\cat{M}} M :=
\mrm{H}^{q} \mrm{R} \Gamma_{\cat{M}} M$.
\end{dfn}

For the purposes of this paper it will be useful to introduce a
notion of flasque modules. Recall that a sheaf $\mcal{M}$ on a
topological space $X$ is called flasque (or flabby) 
if for any two open subsets $V \subset U$ the restriction map
$\Gamma(U, \mcal{M}) \to \Gamma(V, \mcal{M})$
is surjective. It follows that for any closed subset $Z$ the 
cohomology sheaves $\mcal{H}^{q}_{Z} \mcal{M}$, $q > 0$, are zero.
The following definition is somewhat ad
hoc, but we try to justify it in the subsequent examples.

\begin{dfn} \label{dfn1.2}
Suppose $\msf{M} \subset \cat{Mod} A$ is a localizing subcategory. A
module $M$ is called {\em $\msf{M}$-flasque} if
$\mrm{H}^{q}_{\cat{M}} M = 0$ for all $q > 0$.
\end{dfn}

\begin{exa}
Suppose $\msf{M}$ is stable. For any $M \in \msf{M}$ the minimal
injective resolution $M \to I^{0} \to I^{1} \to \cdots$ is in
$\cat{C}^{+}(\msf{M})$, and hence $M$ is $\msf{M}$-flasque.
\end{exa}

\begin{exa} \label{exa1.6}
Suppose $S \subset A$ is a left denominator set (Example
\ref{exa1.3}), and assume that the localizing subcategory
$\msf{M}_{S}$ is stable. Then a module $M$ is
$\msf{M}_{S}$-flasque iff the canonical homomorphism
$M \to A_{S} \otimes_{A} M$ is surjective. To see
why this is true, first observe that for an injective module $I$
the module $\Gamma_{\msf{M}_{S}} I$ is also injective (because of
stability). So $I / \Gamma_{\msf{M}_{S}} I$ is an injective
$S$-torsion-free module. Since $A_{S}$ is a flat right $A$-module
\cite[Proposition II.3.5]{Ste}, we get from
\cite[Proposition V.2.11]{Ste} that
$A_{S} \otimes_{A} I \cong
A_{S} \otimes_{A} (I / \Gamma_{\msf{M}_{S}} I) \cong
I / \Gamma_{\msf{M}_{S}} I$.
Hence there is an exact sequence
$0 \to \Gamma_{\msf{M}_{S}} I \to I \to
A_{S} \otimes_{A} I \to 0$.
Using an injective resolution of $M$ we deduce that
$\mrm{H}^{q}_{\msf{M}_{S}} M = 0$ for $q \geq 2$ and the sequence
$M \to A_{S} \otimes_{A} M \to \mrm{H}^{1}_{\msf{M}_{S}} M \to 0$
is exact.
\end{exa}

A module $M$ is called {\em finitely resolved} if it has a free 
resolution
\[ \cdots \to A^{r_{2}} \to A^{r_{1}} \to A^{r_{0}} \to M \to 0 \]
where all the $r_{i} < \infty$. 

\begin{dfn}
A localizing subcategory $\msf{M}$ is called {\em locally finitely 
resolved} if there is a cofinal inverse system 
$\{ \mfrak{a}_{i} \}$ in the filter $\mfrak{F}_{\msf{M}}$ 
consisting of finitely resolved left ideals.
\end{dfn}

If $A$ is left noetherian then any localizing subcategory $\msf{M}$ 
is automatically locally finitely resolved. But the next examples
shows this is a more general phenomenon.

\begin{exa}
Let $A$ be any ring and $a$ a central regular element. Define
$\mfrak{m} := (a)$. Then the localizing subcategory 
$\msf{M}_{\mfrak{m}}$  (cf.\ Example \ref{exa1.1}) is locally 
finitely resolved. This generalizes to a regular sequence 
$a_{1}, \ldots, a_{n}$ of normalizing elements. 
\end{exa}

\begin{exa}
Let $\mbb{K}$ be a commutative ring and  
$A := \mbb{K} \bra{x_{1}, \ldots, x_{n}}$,
a free associative algebra. Let 
$\mfrak{m} := (x_{1}, \ldots, x_{n})$ be the augmentation ideal. 
Then $\msf{M}_{\mfrak{m}}$ is locally finitely resolved. 
\end{exa}

\begin{exa}
Suppose $A$ is a connected graded algebra over some field 
$\mbb{K}$, and $\mfrak{m}$ is the augmentation ideal. If $A$ is 
Ext-finite in the sense of \cite{VdB1}, i.e.\ every 
$\opn{Ext}^{q}_{A}(\mbb{K}, \mbb{K})$ is finite as 
$\mbb{K}$-module, then the localizing subcategory 
$\msf{M}_{\mfrak{m}}$ is locally finitely resolved. 
This is Van den Bergh's original setup in \cite{VdB1}.
\end{exa}

\begin{prop} \label{prop1.6}
Suppose $\msf{M}$ is locally finitely resolved. Then
$\mrm{H}^{q}_{\msf{M}}$ commutes with direct limits. Therefore
the direct limit of $\msf{M}$-flasque modules is 
$\msf{M}$-flasque.
\end{prop}

\begin{proof}
Let $\{ \mfrak{a}_{i} \}$ be a cofinal inverse system in the
filter $\mfrak{F}_{\msf{M}}$ with all the left ideals
$\mfrak{a}_{i}$ finitely resolved. Say
$M = \underset{j \to}{\lim}\, M_{j}$ for some direct system 
$\{ M_{j} \}$ of $A$-modules. 
Since the left module $A / \mfrak{a}_{i}$ is finitely resolved we get
\[ \opn{Ext}^{q}_{A}(A / \mfrak{a}_{i}, M) \cong
\underset{j \to}{\lim} \opn{Ext}^{q}_{A}(A / \mfrak{a}_{i}, M_{j}) 
. \]
Hence for any $q$ we have
\[ \begin{aligned}
\mrm{H}^{q}_{\msf{M}} M & \cong 
\underset{i \to}{\lim} \opn{Ext}^{q}_{A}(A / \mfrak{a}_{i}, M) \\
& \cong \underset{i \to}{\lim} \, \underset{j \to}{\lim} 
\opn{Ext}^{q}_{A}(A / \mfrak{a}_{i}, M_{j}) \\
& \cong \underset{j \to}{\lim} \, \underset{i \to}{\lim} 
\opn{Ext}^{q}_{A}(A / \mfrak{a}_{i}, M_{j}) \\
& \cong \underset{j \to}{\lim} \, \mrm{H}^{q}_{\msf{M}} M_{j} .
\end{aligned} \]
\end{proof}

Since any injective module is $\msf{M}$-flasque it follows that
there are enough $\msf{M}$-flasque modules: any module embeds into
an $\msf{M}$-flasque one. Hence for any
$M \in \cat{D}^{+}(\cat{Mod} A)$ there is a quasi-isomorphism
$M \to I$ in $\cat{C}^{+}(\cat{Mod} A)$ with each $I^{q}$ an
$\msf{M}$-flasque module. We call such a quasi-isomorphism an
{\em $\msf{M}$-flasque resolution of $M$ in
$\cat{C}^{+}(\cat{Mod} A)$}.

\begin{prop} \label{prop1.4}
Let $M \in \cat{D}^{+}(\cat{Mod} A)$ and
$M \to I$ an $\msf{M}$-flasque resolution in
$\cat{C}^{+}(\cat{Mod} A)$. Then the canonical morphism
$\Gamma_{\cat{M}} I \to \mrm{R} \Gamma_{\cat{M}} I$
is an isomorphism, and hence
$\mrm{R} \Gamma_{\cat{M}} M \cong \Gamma_{\cat{M}} I$.
\end{prop}

\begin{proof}
If $J \in \cat{C}^{+}(\cat{Mod} A)$ is an acyclic
complex of $\msf{M}$-flasque modules then
$\Gamma_{\cat{M}} J$ is also acyclic. Now use
\cite[Theorem I.5.1]{RD}.
\end{proof}

Thus we can compute $\mrm{R} \Gamma_{\cat{M}}$ using $\cat{M}$-flasque
resolutions.

Let $\mbb{K}$ be a commutative base ring and let $A$ and $B$ be
associative unital $\mbb{K}$-algebras. We denote by $B^{\mrm{op}}$ the
opposite algebra, and
$A \otimes B^{\mrm{op}} := A \otimes_{\mbb{K}} B^{\mrm{op}}$. Thus
an $(A \otimes B^{\mrm{op}})$-module $M$ is, in conventional notation,
a $\mbb{K}$-central $A$-$B$-bimodule ${}_{A}M_{B}$. When $A = B$
then $A^{\mrm{e}} := A \otimes A^{\mrm{op}}$ is the enveloping
algebra.

\begin{prop} \label{prop1.3}
Let $A$ and $B$ be $\mbb{K}$-algebras with $B$ flat over 
$\mbb{K}$. Then there is a derived functor
\[ \mrm{R} \Gamma_{\cat{M}} :
\cat{D}^{+}(\cat{Mod} A \otimes B^{\mrm{op}})
\to \cat{D}^{+}(\cat{Mod} A \otimes B^{\mrm{op}})  \]
commuting with the forgetful functor
$\cat{D}^{+}(\cat{Mod} A \otimes B^{\mrm{op}}) \to
\cat{D}^{+}(\cat{Mod} A)$. 
In particular an $A \otimes B^{\mrm{op}}$-module $M$ is 
$\cat{M}$-flasque iff it is $\cat{M}$-flasque as $A$-module.
\end{prop}

\begin{proof}
Since $B^{\mrm{op}}$ is flat
over $\mbb{K}$, any injective $A \otimes B^{\mrm{op}}$-module is
also an injective $A$-module.
\end{proof}

The next theorem is inspired by \cite[Theorem 4.8]{VdB1}. We shall
use it in our discussion of Cousin complexes in the next section.

\begin{thm} \label{thm1.1}
Let $A$ and $B$ be flat $\mbb{K}$-algebras and
let $\cat{M} \subset \cat{Mod} A$ and
$\cat{N} \subset \cat{Mod} B^{\mrm{op}}$
be stable, locally finitely resolved, localizing subcategories. Suppose
$M \in \cat{D}^{+}(\cat{Mod} A \otimes B^{\mrm{op}})$
satisfies
$\mrm{H}^{q}_{\cat{M}} M \in \cat{N}$
and
$\mrm{H}^{q}_{\cat{N}} M \in \cat{M}$
for all integers $q$. Then there is a functorial isomorphism
\[ \mrm{R} \Gamma_{\cat{M}} M \cong
\mrm{R} \Gamma_{\cat{N}} M \text{ in }
\cat{D}(\cat{Mod} A \otimes B^{\mrm{op}}) . \]
\end{thm}

We precede the proof by three lemmas.

\begin{lem} \label{lem1.2}
In the situation of the theorem, but without the stability assumption,
if $I$ is an injective $A \otimes B^{\mrm{op}}$-module, then
$\Gamma_{\cat{M}} I$ is an $\msf{N}$-flasque $B^{\mrm{op}}$-module.
\end{lem}

\begin{proof}
Let $\mfrak{F}_{\cat{M}}$ be the filter of left ideals associated
with $\cat{M}$. Then
\[ \Gamma_{\cat{M}} I = \varinjlim_{\mfrak{a} \in
\mfrak{F}_{\cat{M}}} \opn{Hom}_{A}(A / \mfrak{a}, I) . \]
Because $B^{\mrm{op}}$ is flat over $\mbb{K}$ each
$\opn{Hom}_{A}(A / \mfrak{a}, I)$ is an injective
$B^{\mrm{op}}$-module, and hence it is $\msf{N}$-flasque. By 
Proposition \ref{prop1.6} the direct limit of $\msf{N}$-flasque 
modules is $\msf{N}$-flasque.
\end{proof}

\begin{lem} \label{lem1.3}
In the situation of the theorem, but without the stability assumption,
there is a functorial isomorphism
\[ \mrm{R} \Gamma_{\cat{N}} \mrm{R} \Gamma_{\cat{M}} M
\cong \mrm{R}(\Gamma_{\cat{N}} \Gamma_{\cat{M}}) M
\text{ in } \cat{D}^{+}(\cat{Mod} A \otimes B^{\mrm{op}}) . \]
\end{lem}

\begin{proof}
Take an injective resolution $M \to I$ in
$\cat{C}^{+}(\cat{Mod} A \otimes B^{\mrm{op}})$.
Then 
$\mrm{R} (\Gamma_{\cat{N}} \Gamma_{\cat{M}}) M$ \linebreak
$= \Gamma_{\cat{N}} \Gamma_{\cat{M}} I$, 
and
$\mrm{R} \Gamma_{\cat{M}} M = \Gamma_{\cat{M}} I$.
According to Lemma \ref{lem1.2},
$\Gamma_{\cat{M}} I$ is a complex of $\msf{N}$-flasque
$B^{\mrm{op}}$ modules, so by Proposition \ref{prop1.4},
$\mrm{R} \Gamma_{\cat{N}} \Gamma_{\cat{M}} I =
\Gamma_{\cat{N}} \Gamma_{\cat{M}} I$. 
\end{proof}

\begin{lem} \label{lem1.4}
In the situation of the theorem, let
$N \in \cat{D}^{+}_{\cat{N}}(\cat{Mod} B^\mrm{op})$.
Then the natural morphism
$\mrm{R} \Gamma_{\cat{N}} N \to N$
in $\cat{D}(\cat{Mod} B^\mrm{op})$ is an isomorphism.
\end{lem}

\begin{proof}
It suffices (by ``way-out'' reasons, see
\cite[Proposition 7.1(iv)]{RD})
to consider a single module $N \in \cat{N}$.
But by the stability assumption, the minimal injective resolution
$N \to I^{0} \to I^{1} \to \cdots$ is in $\cat{N}$, so
$I = \Gamma_{\cat{N}} I$.
\end{proof}

\begin{proof}[Proof of Theorem \tup{\ref{thm1.1}}]
For any bimodule $N$ write
\[ \Gamma_{\cat{M} \cap \cat{N}} N :=
\Gamma_{\cat{N}} \Gamma_{\cat{M}} N =
\Gamma_{\cat{N}} N \cap \Gamma_{\cat{M}} N =
\Gamma_{\cat{M}} \Gamma_{\cat{N}} N \subset N, \]
It suffices by symmetry to prove that
$\mrm{R} \Gamma_{\cat{M}} M \cong
\mrm{R} \Gamma_{\cat{M} \cap \cat{N}} M$.
Since
$\mrm{R} \Gamma_{\cat{M}} M \in
\cat{D}^{+}_{\cat{N}}(\cat{Mod} B^\mrm{op})$,
Lemma \ref{lem1.4} says that there is a functorial isomorphism
$\mrm{R} \Gamma_{\cat{M}} M \cong
\mrm{R} \Gamma_{\cat{N}} \mrm{R} \Gamma_{\cat{M}} M$.
Finally by Lemma \ref{lem1.3} there is a functorial isomorphism
$\mrm{R} \Gamma_{\cat{M} \cap \cat{N}} M \cong
\mrm{R} \Gamma_{\cat{N}} \mrm{R} \Gamma_{\cat{M}} M$.
\end{proof}

\begin{exa} \label{exa1.5}
Theorem \ref{thm1.1} does not hold in general without the stability
assumption. Here is a counterexample. Take $\mbb{K} = \mbb{C}$ and
$A = B = \mrm{U}(\mfrak{sl}_2)$, the universal enveloping algebra
of the Lie algebra $\mfrak{sl}_2$. Let 
$\cat{M}_{0} := \cat{M}_{0}(\opn{GKdim})$ 
be the full subcategory of $\cat{Mod} A$ consisting of 
modules of Gelfand-Kirillov dimension $0$ (unions of $A$-modules that 
are finite over $\mbb{K}$), and let $\cat{N} := \cat{Mod} B^\mrm{op}$. 
Then all hypotheses of Theorem \ref{thm1.1} hold,
except that $\cat{M}_{0}$ is not stable.
If we take $M \in \cat{M}_{0}$ to be the simple module with 
$\opn{rank}_{\mbb{K}} M = 1$, then it follows from
\cite[Proposition 7.5]{AjSZ} that
$\mrm{H}^{3}_{\cat{M}_{0}} M \cong M$, but of course
$\mrm{H}^{3}_{\cat{N}} M = 0$.
\end{exa}

\section{Cousin Functors}

Cousin complexes in commutative algebraic geometry were introduced 
by Groth\-en\-dieck in \cite{RD}. Several people (including Lipman, 
private communication) had suggested extending the construction to 
more general settings. The noncommutative version below 
already appeared in \cite{Ye2}, but the powerful Theorem 
\ref{thm2.1} is new. 

Suppose $A$ is a ring, and we are given an increasing filtration
\[ \cdots \subset  \cat{M}_{q - 1} \subset \cat{M}_{q} \subset
\cat{M}_{q + 1} \subset \cdots \subset \cat{Mod} A \]
by localizing subcategories, indexed by $\mbb{Z}$. We shall sometimes
write $\cat{M}^{q} := \cat{M}_{-q}$, so that
$\{ \cat{M}^{p} \}_{p \in \mbb{Z}}$ is a
decreasing filtration. This is to conform to the convention that
decreasing filtrations go with cochain complexes.
We say the filtration
$\cat{M} = \{ \cat{M}_{q} \} = \{ \cat{M}^{p} \}$
is  {\em bounded} if there are $q_{0} \leq q_{1}$ such that
$\cat{M}_{q_{0} - 1} = 0$ and $\cat{M}_{q_{1}} = \cat{Mod} A$.

\begin{exa} \label{exa2.1}
Suppose $\opn{dim}$ is an exact dimension function that is bounded, in
the sense that there are integers $q_{0} \leq q_{1}$ such that for
any nonzero module $M$, $q_{0} \leq \opn{dim} M \leq q_{1}$.
Define $\cat{M}_{q}(\opn{dim})$ as in Example \ref{exa1.2}. Then 
$\cat{M} = \{ \cat{M}_{q}(\opn{dim}) \}$ is a bounded filtration of 
$\cat{Mod} A$ by localizing subcategories.
Conversely, given a bounded filtration 
$\cat{M} = \{ \cat{M}_{q} \}$ by localizing subcategories,
we can define 
$\opn{dim} M := \opn{inf} \{ q \mid M \in \cat{M}_{q} \}$, 
and this will be a bounded exact dimension function.
\end{exa}

\begin{exa}  \label{exa2.2}
Specializing the previous example, let $A$ be a finitely generated
commutative algebra over a field $\mbb{K}$ (or more generally $A$ is a
catenary commutative noetherian ring of Krull dimension
$\opn{Kdim} A < \infty$) and $X := \opn{Spec} A$.
Taking $\opn{dim} = \opn{Kdim}$,
$\mcal{Z}_{q} = \{ Z \subset X \text{ closed, } \opn{dim} Z \leq q \}$
and $\mcal{Z}^{q} := \mcal{Z}_{-q}$ we get
$\Gamma_{\cat{M}^{p}} = \Gamma_{\mcal{Z}^{p}}$.
This is the kind of filtration by codimension (coniveau) used in
\cite[Chapter IV]{RD}. The bounds are $q_{1} = \opn{dim} X$ and
$q_{0} = 0$.
\end{exa}

Suppose $\msf{M} = \{ \msf{M}_{q} \}$ is a collection
of localizing subcategories of $\cat{Mod} A$. We call
a module $M$ {\em $\msf{M}$-flasque} if it is
$\msf{M}_{q}$-flasque for all $q$ (Definition \ref{dfn1.2}).

For a module $M$ and $d \geq 0$ we write
$\Gamma_{\cat{M}^{p} /  \cat{M}^{p + d}} M :=
\Gamma_{\cat{M}^{p}} M / \Gamma_{\cat{M}^{p + d}} M$.

\begin{lem} \label{lem2.1}
There is a right derived functor
$\mrm{R} \Gamma_{\cat{M}^{p} / \cat{M}^{p + d}}$
that fits into a functorial triangle for
$M \in \msf{D}^{+}(\cat{Mod} A)$:
\[ \mrm{R} \Gamma_{\cat{M}^{p + d}} M \to
\mrm{R} \Gamma_{\cat{M}^{p}} M \to
\mrm{R} \Gamma_{\cat{M}^{p} /  \cat{M}^{p + d}} M \to
\mrm{R} \Gamma_{\cat{M}^{p + d}} M[1] . \]
If $M \to I$ is a flasque resolution then 
$\mrm{R} \Gamma_{\cat{M}^{p} /  \cat{M}^{p + d}} M =
\Gamma_{\cat{M}^{p} /  \cat{M}^{p + d}} I$.
\end{lem}

\begin{proof}
If $I \in \msf{D}^{+}(\cat{Mod} A)$ is an acyclic complex of $
\msf{M}$-flasque modules then from the exact sequence of complexes
\[ 0 \to \Gamma_{\cat{M}^{p + d}} I \to
\Gamma_{\cat{M}^{p}} I \to
\Gamma_{\cat{M}^{p} /  \cat{M}^{p + d}} I \to 0 \]
we see that 
$\Gamma_{\cat{M}^{p} /  \cat{M}^{p + d}} I$
is also acyclic. Thus we can define 
$\mrm{R} \Gamma_{\cat{M}^{p} /  \cat{M}^{p + d}} M :=
\Gamma_{\cat{M}^{p} /  \cat{M}^{p + d}} I$
when $M \in \msf{D}^{+}(\cat{Mod} A)$ and 
$M \to I$ is a flasque resolution (cf.\ proof of Proposition 
\ref{prop1.4}).  
\end{proof}

We set
$\mrm{H}^{q}_{\cat{M}^{p} /  \cat{M}^{p + d}} M :=
\mrm{H}^{q} \mrm{R} \Gamma_{\cat{M}^{p} /  \cat{M}^{p + d}} M$.

\begin{prop} \label{prop2.1}
Let $\cat{M} = \{ \cat{M}^{p} \}$ be a bounded filtration of
$\cat{Mod} A$ by localizing subcategories. Then for every
$M \in \msf{D}^{+}(\cat{Mod} A)$ there is a convergent spectral
sequence
\[  E_{1}^{p, q} =
\mrm{H}^{p + q}_{\cat{M}^{p} /  \cat{M}^{p + 1}} M
\Rightarrow \mrm{H}^{p + q} M , \]
functorial in $M$.
\end{prop}

\begin{proof}
Pick an $\msf{M}$-flasque resolution $M \to I$ in
$\msf{C}^{+}(\cat{Mod} A)$.
The decreasing filtration $\{ \Gamma_{\cat{M}^{p}} I \}$ is
bounded in the sense of \cite[Sections XI.3 and XI.8]{Mac}, i.e.\
$\Gamma_{\cat{M}^{p_{0}}} I = I$ and
$\Gamma_{\cat{M}^{p_{1}}} I = 0$ for some $p_{0} \leq p_{1}$.
Hence by \cite[Theorem XI.3.1]{Mac} we get a convergent spectral
sequence
\[ E_{1}^{p, q} =
\mrm{H}^{p + q} \Gamma_{\cat{M}^{p} /  \cat{M}^{p + 1}} I
\Rightarrow \mrm{H}^{p + q} I . \]
Now
$\mrm{H}^{p + q} \Gamma_{\cat{M}^{p} / \cat{M}^{p + 1}} I =
\mrm{H}^{p + q}_{\cat{M}^{p} / \cat{M}^{p + 1}} M$
and
$\mrm{H}^{p + q} I = \mrm{H}^{p + q} M$.
If $I_{1} \to I_{2}$ is a homomorphism between $\msf{M}$-flasque
complexes then there is a map between the two spectral sequences;
and if $I_{1} \to I_{2}$ is a quasi-isomorphism then the two spectral
sequences are isomorphic.
\end{proof}

\begin{dfn}(Grothendieck, \cite{RD}) \label{dfn2.1}
Given a bounded filtration $\cat{M} = \{ \cat{M}^{p} \}$ of
$\cat{Mod} A$ and a complex
$M \in \cat{D}^{+}(\cat{Mod} A)$ define the complex
$\mrm{E}_{\cat{M}} M$ as follows. For any $p$
\[ (\mrm{E}_{\cat{M}} M)^{p} := E_{1}^{p, 0} =
\mrm{H}^{p}_{\cat{M}^{p} /  \cat{M}^{p + 1}} M  \]
in the spectral sequence above, and the coboundary operator is
\[ \mrm{d}^{p, 0}_{1} : (\mrm{E}_{\cat{M}} M)^{p} = E^{p, 0}_{1}
\to (\mrm{E}_{\cat{M}} M)^{p + 1} = E^{p + 1, 0}_{1} . \]
Thus $\mrm{E}_{\cat{M}} M$ is the row $q = 0$ in the $E_{1}$
page of the spectral sequence. We obtain an additive functor
\[ \mrm{E}_{\cat{M}} :
\cat{D}^{+}(\cat{Mod} A) \to \cat{C}^{+}(\cat{Mod} A) \]
called the {\em Cousin functor}.
\end{dfn}

Unlike the commutative situation, here the complex
$\mrm{E}_{\msf{M}} M$ can behave quite oddly -- see below.

\begin{dfn} \label{dfn2.3}
Given an exact dimension function $\opn{dim}$ on $\cat{Mod} A$
we say an $A$-module $M$ is 
{\em $\opn{dim}$-pure of dimension $q$}
if $\opn{dim} M' = \opn{dim} M = q$ for all nonzero submodules
$M' \subset M$.
\end{dfn}

\begin{rem} \label{rem2.1}
In the commutative case (see Example \ref{exa2.2}) let
$M \in \msf{D}^{+}(\cat{Mod} A)$ and let
$\mcal{M} := \mcal{O}_{X} \otimes_{A} M$ denote the corresponding
complex of quasi-coherent sheaves on $X$. Then for any $p, q$ one has
\[ \mrm{H}^{p + q}_{\cat{M}^{p} / \cat{M}^{p + 1}} M \cong
\Gamma(X, \mcal{H}^{p + q}_{\mcal{Z}^{p} / \mcal{Z}^{p + 1}}
\mcal{M}) \cong
\bigoplus_{x} \mrm{H}_{x}^{p + q} \mcal{M} \]
where $x$ runs over the points in $X$ of
$\opn{dim} \overline{\{ x \}} = -p$
and $\mrm{H}_{x}^{p + q}$ is local cohomology. In the language of
\cite{RD}, the sheaf
$\mcal{H}^{p + q}_{\mcal{Z}^{p} / \mcal{Z}^{p + 1}} \mcal{M}$
lies on the $\mcal{Z}^{p} / \mcal{Z}^{p + 1}$-skeleton of $X$. In
particular this means the $A$-module
$(\mrm{E}_{\msf{M}} M)^{-q} =
\mrm{H}^{-q}_{\cat{M}_{q} / \cat{M}_{q - 1}} M$
is $\msf{M}$-flasque and $\opn{Kdim}$-pure of dimension $q$.
Note that this implies
$\mrm{E}_{\msf{M}} \mrm{E}_{\msf{M}} M = \mrm{E}_{\msf{M}} M$.
A complex $N$ such that each $N^{-q}$ is $\msf{M}$-flasque and
$\opn{Kdim}$-pure of dimension $q$ is called a ``Cousin complex'' in
\cite[Section IV.3]{RD}. However for a noncommutative ring $A$ the
complex $\mrm{E}_{\msf{M}} M$ will seldom be a Cousin complex in
this sense; cf.\ the next example.
\end{rem}

\begin{exa} \label{exa2.4}
Consider $\mbb{K} = \mbb{C}$ and $A = \mrm{U}(\mfrak{sl}_2)$ as
in Example \ref{exa1.5}. Let
$\msf{M} = \{ \msf{M}_{q}(\opn{GKdim}) \}$ be the filtration
by Gelfand-Kirillov dimension and $M \in \msf{M}_{0}$ 
the simple $A$-module. Then
$(\mrm{E}_{\msf{M}} M)^{0} =
\mrm{H}^{0}_{\msf{M}_{0}} M \cong M$.
Since $\mrm{H}^{3}_{\msf{M}_{0}} M \neq 0$ we see that
$(\mrm{E}_{\msf{M}} M)^{0}$
is not $\msf{M}$-flasque.
\end{exa}

\begin{prop} \label{prop2.6}
Suppose $A \to B$ is a homomorphism of rings,
$\msf{M}(A) = \{ \msf{M}_{q}(A) \}$ and
$\msf{M}(B) = \{ \msf{M}_{q}(B) \}$ are bounded filtrations of
$\cat{Mod} A$ and $\cat{Mod} B$ respectively by localizing
subcategories, with Cousin functors $\mrm{E}_{\msf{M}(A)}$ and
$\mrm{E}_{\msf{M}(B)}$, satisfying:
\begin{enumerate}
\rmitem{i} For any $M \in \cat{Mod} B$ and any $q$,
$\Gamma_{\cat{M}_{q}(B)} M = \Gamma_{\cat{M}_{q}(A)} M$.
\rmitem{ii} If $I \in \cat{Mod} B$ is injective then
it is $\msf{M}(A)$-flasque.
\end{enumerate}
Then there is an isomorphism
$\mrm{E}_{\msf{M}(B)} M \cong \mrm{E}_{\msf{M}(A)} M$,
functorial in $M \in \cat{D}^{+}(\cat{Mod} B)$.
\end{prop}

\begin{proof}
Choose an injective resolution $M \to I$ in
$\cat{C}^{+}(\cat{Mod} B)$. Then $M \to I$ is an $\msf{M}(A)$-flasque
resolution in $\cat{C}^{+}(\cat{Mod} A)$, the filtered
complexes $\Gamma_{\cat{M}_{q}(B)} I$ and
$\Gamma_{\cat{M}_{q}(A)} I$ coincide, and the spectral sequence
defines both $\mrm{E}_{\msf{M}(B)} M$ and $\mrm{E}_{\msf{M}(A)} M$.
\end{proof}

Now let $\mbb{K}$ be a commutative base ring and as before
$\otimes = \otimes_{\mbb{K}}$.

\begin{cor} \label{cor2.5}
Let $A$ be a $\mbb{K}$-algebra, $\msf{M} = \{ \msf{M}^{p} \}$
a bounded filtration of $\cat{Mod} A$ by localizing subcategories,
and $B$ a flat $\mbb{K}$-algebra. 
Suppose $M \in \cat{D}^{+}(\cat{Mod} A \otimes B^{\mrm{op}})$.
Then the Cousin functor $\mrm{E}_{\cat{M}} M$ commutes with the 
forgetful functor
$\cat{Mod} A \otimes B^{\mrm{op}} \to \cat{Mod} A$.
\end{cor}

\begin{proof}
Write $\msf{M}(A)$ and $\msf{M}(A \otimes B^{\mrm{op}})$ for the 
filtrations of $\cat{Mod} A$  and $\cat{Mod} A \otimes B^{\mrm{op}}$
respectively, and apply Proposition \ref{prop2.6}.
\end{proof}

For the next theorem it will be important to distinguish between
morphisms in $\cat{D}(\cat{Mod} A)$ and $\cat{C}(\cat{Mod} A)$,
so let
$\mrm{Q} : \cat{C}(\cat{Mod} A) \to \cat{D}(\cat{Mod} A)$
be the localization functor (identity on objects).

\begin{thm} \label{thm2.1}
Let $A$ be a ring, $\cat{M} = \{ \cat{M}^{p} \}$ a bounded
filtration of $\cat{Mod} A$ by localizing subcategories and
$\mrm{E}_{\cat{M}}$ the associated Cousin functor.
Let $M \in \cat{D}^{\mrm{b}}(\cat{Mod} A)$ a complex satisfying
\begin{enumerate}
\rmitem{$\ast$}
$\mrm{H}^{p + q}_{\cat{M}^{p} /  \cat{M}^{p + 1}} M = 0$
for all $q \neq 0$ and all $p$.
\end{enumerate}
Then there is an isomorphism
$M \cong \mrm{Q} \mrm{E}_{\cat{M}} M$ in $\cat{D}(\cat{Mod} A)$.
\end{thm}

\begin{proof}
This is really the implication (ii) $\Rightarrow$ (iii) in
\cite[Proposition IV.3.1]{RD}. We shall explain the minor modification
needed in the proof to make it apply to our situation. Also we shall
sketch the main ideas of the proof using our notation, so the
interested reader can find it easier to consult the rather lengthy
proof in \cite{RD}.

The result in \cite{RD} refers to the abelian category
$\cat{Ab} X$ of sheaves of abelian groups on a topological space $X$.
The space $X$ has a filtration $\{ Z^{p} \}$ by closed subsets,
inducing a filtration
$\cat{M} = \{ \cat{M}^{p} \}$ of $\cat{Ab} X$, with
$\Gamma_{\cat{M}^{p}} = \Gamma_{Z^{p}}$.
With this notation the proof involves homological algebra only,
hence it applies to $\cat{Mod} A$ as well.

Here is the sketch. Let us abbreviate
$\mrm{E} := \mrm{E}_{\msf{M}}$. Define
\[ \tau_{\geq p} \mrm{E} M :=
(\cdots \to 0 \to (\mrm{E} M)^{p} \to (\mrm{E} M)^{p + 1}
\to \cdots) \]
to be the truncation in $\msf{C}(\cat{Mod} A)$. One shows by
descending induction on $p$ that there are (noncanonical)
isomorphisms
\begin{equation} \label{eqn2.1}
\phi_{p}: \mrm{R} \Gamma_{\cat{M}^{p}} M \iso
\mrm{Q} \tau_{\geq p} \mrm{E} M \text{ in } \msf{D}(\cat{Mod} A)
\end{equation}
such that the diagrams
\begin{equation} \label{eqn2.2}
\choosegraphics{\UseTips \xymatrix{
{\mrm{H}}^{p}_{\cat{M}^{p} / \cat{M}^{p + 1}} M \ar[d]^{=} &
{\mrm{H}}^{p}_{\cat{M}^{p}} M \ar[l] \ar[d]^{\mrm{H}^{p}(\phi_{p})} \\
(\mrm{E} M)^{p}  &
{\mrm{H}}^{p} \tau_{\geq p} \mrm{E} M \ar[l]
}}
\end{equation}
commute. The horizontal arrows are the canonical ones.

The starting point is that for large enough $p = p_{\mrm{big}}$,
$\mrm{H}^{p} M = \mrm{H}^{p} \mrm{E} M = 0$. For such $p$ one
shows that $\mrm{H}^{q}_{\cat{M}^{p}} M = 0$ if $q \neq p$. Hence
there is an isomorphism $\mrm{R} \Gamma_{\cat{M}^{p}} M \cong
\mrm{Q} (\mrm{H}^{p}_{\cat{M}^{p}} M)[-p]$. Also one shows that
$\mrm{H}^{q} \tau_{\geq p} \mrm{E} M = 0$ if $q \neq p$, so that
$\mrm{Q} \tau_{\geq p} \mrm{E} M \cong \mrm{Q} (\mrm{H}^{p}
\tau_{\geq p} \mrm{E} M)[-p]$. Since the modules
$\mrm{H}^{p}_{\cat{M}^{p}} M$ and $\mrm{H}^{p} \tau_{\geq p}
\mrm{E} M$ are canonically isomorphic in this case, we get an
isomorphism (\ref{eqn2.1}) for $p = p_{\mrm{big}}$.

In the inductive step, depicted in Figure \ref{fig1}, 
we have two canonical triangles (in which the morphisms 
$\alpha_{p - 1}$ and
$\mrm{d}^{p - 1}$ have degree $+1$), canonical isomorphisms
$\psi_{p}$ and $\psi_{p - 1}$ (arising from the assumption
($\ast$)) and an isomorphism $\phi_{p}$ that's already been
constructed. The square on the left commutes because diagram
(\ref{eqn2.2}) commutes.

\begin{figure}
\[ \choosegraphics{\UseTips \xymatrix{
{\mrm{R}} \Gamma_{\cat{M}^{p} / \cat{M}^{p + 1}} M
\ar[dd]^{\psi_{p}}
& & {\mrm{R}} \Gamma_{\cat{M}^{p - 1} / \cat{M}^{p}} M
\ar[dl]_{\alpha_{p - 1}} \ar '[d] [dd]^{\psi_{p - 1}} \\
& {\mrm{R}} \Gamma_{\cat{M}^{p}} M \ar[rr] \ar[dd]^{\phi_{p}}
\ar[ul]_{\beta_{p}} & &
{\mrm{R}} \Gamma_{\cat{M}^{p - 1}} M \ar[ul]_{\beta_{p - 1}}
\ar@{-->}[dd]^{\phi_{p - 1}} \\
{\mrm{Q}} (\mrm{E} M)^{p}[-p]
& & {\mrm{Q}} (\mrm{E} M)^{p - 1}[-p + 1] \ar[dl]_{\mrm{d}^{p - 1}} \\
& {\mrm{Q}} \tau_{\geq p} \mrm{E} M \ar[rr]  \ar[ul]
& & {\mrm{Q}} \tau_{\geq p - 1} \mrm{E} M  \ar[ul]
}} \]
\caption{} 
\label{fig1}
\end{figure}

By definition  of the Cousin complex it follows that
$\mrm{d}^{p - 1} = \mrm{H}^{p - 1}(\beta_{p} \alpha_{p - 1})$.
Since
$\mrm{H}^{p} \tau_{\geq p} \mrm{E} M \subset (\mrm{E} M)^{p}$,
diagram (\ref{eqn2.2}) implies that
\[ \mrm{H}^{p}(\phi_{p} \alpha_{p - 1}) =
\mrm{H}^{p - 1}(\mrm{d}^{p - 1} \psi_{p - 1}):
\mrm{H}^{p - 1}_{\cat{M}^{p - 1} / \cat{M}^{p}} M \to
\mrm{H}^{p} \tau_{\geq p} \mrm{E} M . \]
Therefore $\phi_{p} \alpha_{p - 1} = d^{p - 1} \psi_{p - 1}$,
and hence, by the axioms of triangulated categories, 
there is an isomorphism $\phi_{p - 1}$ making diagram in Figure 
\ref{fig1} commutative. Note that diagram (\ref{eqn2.2})
for $p - 1$ commutes too so the induction continues.
\end{proof}

\begin{rem}
In \cite{RD} a complex satisfying condition ($\ast$) of the theorem is
called a Cohen-Macaulay complex w.r.t.\ the filtration. And indeed
in the commutative case (Example \ref{exa2.2}), for an $A$-module
$M$ of $\opn{Kdim} M = d$, $M$ is a Cohen-Macaulay module iff the
complex $M[d]$ satisfies ($\ast$); cf. \cite[page 239]{RD}.
For a noncommutative ring $A$ these notions diverge.
\end{rem}

\begin{cor} \label{cor2.2}
Let $A$ be a $\mbb{K}$-algebra, $\msf{M} = \{ \msf{M}^{p} \}$
a bounded filtration of $\cat{Mod} A$ by localizing subcategories,
and $B$ a flat $\mbb{K}$-algebra.
Suppose $M \in \cat{D}^{+}(\cat{Mod} A \otimes B^{\mrm{op}})$
satisfies condition \tup{($\ast$)} of the theorem.
Then there is an isomorphism
$M \cong \mrm{Q} \mrm{E}_{\cat{M}} M$ in
$\cat{D}^{+}(\cat{Mod} A \otimes B^{\mrm{op}})$
commuting with the forgetful functor
$\cat{Mod} A \otimes B^{\mrm{op}} \to \cat{Mod} A$.
\end{cor}

\begin{proof}
Invoke the theorem with $A \otimes B^{\mrm{op}}$
instead of $A$, and use Corollary \ref{cor2.5}.
\end{proof}

We shall also need the next propositions.

\begin{prop} \label{prop2.2}
Suppose $A$ and $B$ are  flat $\mbb{K}$-algebras, and
$\cat{M} = \{ \cat{M}^{p} \}$ and $\cat{N} = \{ \cat{N}^{p} \}$
are bounded filtrations of $\cat{Mod} A$ and
$\cat{Mod} B^{\mrm{op}}$ respectively by stable, locally 
finitely resolved, localizing subcategories. Let
$M \in \cat{D}^{+}(\cat{Mod} A \otimes B^{\mrm{op}})$
be a complex satisfying
$\mrm{H}^{q}_{\cat{M}^{p}} M \in \cat{N}^{p}$
and
$\mrm{H}^{q}_{\cat{N}^{p}} M \in \cat{M}^{p}$
for all $p, q$. Then there is a functorial isomorphism
\[  \mrm{E}_{\cat{M}} M \cong \mrm{E}_{\cat{N}} M
\text{ in } \cat{C}(\cat{Mod} A \otimes B^{\mrm{op}})  . \]
\end{prop}

\begin{proof}
Choose an injective resolution $M \to I$ in
$\cat{C}^{+}(\cat{Mod} A \otimes B^{\mrm{op}})$.
Denote by $\Gamma_{\cat{M}} I$ the filtered complex with
filtration
$\{ \Gamma_{\cat{M}^{p}} I \}_{p \in \mbb{Z}}$,
and by $\Gamma_{\cat{M} \cap \cat{N}} I$ the filtered complex with
filtration
$\{ \Gamma_{\cat{M}^{p}} \Gamma_{\cat{N}^{p}}
I \}_{p \in \mbb{Z}}$.
By the proof of Theorem \ref{thm1.1},
the homomorphism of filtered complexes
$\Gamma_{\cat{M} \cap \cat{N}} I \to \Gamma_{\cat{M}} I$
induces an isomorphism on the $E_{1}$ pages of the spectral of the
sequences from Proposition \ref{prop2.1}. Similarly for
$\Gamma_{\cat{M} \cap \cat{N}} I \to \Gamma_{\cat{N}} I$.
\end{proof}
%

\begin{prop} \label{prop2.7}
Let $\opn{dim}$ be a bounded exact dimension function on 
$\cat{Mod} A$, and let 
$\msf{M} = \{ \msf{M}_{q}(\dim) \}$ be a the induced filtration of 
$\cat{Mod} A$. Suppose the complexes
$M, I \in \cat{C}^{+}(\cat{Mod} A)$ satisfy:
\begin{enumerate}
\rmitem{i} Each module $M^{-q}$ and $I^{-q}$ is $\msf{M}$-flasque
and $\opn{dim}$-pure of dimension $q$.
\rmitem{ii} Each module $I^{-q}$ is injective.
\end{enumerate}
Then
\begin{enumerate}
\item $\mrm{E}_{\msf{M}} \mrm{Q} M = M$ and
$\mrm{E}_{\msf{M}} \mrm{Q} I = I$.
\item The functor $\mrm{E}_{\msf{M}}$ induces an isomorphism
\[ \opn{Hom}_{\cat{D}^{+}(\cat{Mod} A)}(\mrm{Q} M, \mrm{Q} I)
\iso \opn{Hom}_{\cat{C}^{+}(\cat{Mod} A)}(M,I) \]
with inverse induced by $\mrm{Q}$.
\end{enumerate}
\end{prop}

\begin{proof}
1. Clear, since
$\Gamma_{\msf{M}^{p} / \msf{M}^{p + 1}} M^{-p} \cong M^{-p}$
and the same for $I$.

\medskip \noindent 2. Since $I$ is a bounded below complex of
injectives we have an isomorphism
\[ \mrm{H}^{0} \opn{Hom}_{A}(M,I) \cong
\opn{Hom}_{\cat{D}^{+}(\cat{Mod} A)}(\mrm{Q} M, \mrm{Q} I) . \]
The purity implies that $\opn{Hom}_{A}(M,I)^{-1} = 0$
and hence we get an isomorphism
\[ \opn{Hom}_{\cat{C}^{+}(\cat{Mod} A)}(M, I) \iso
\opn{Hom}_{\cat{D}^{+}(\cat{Mod} A)}(\mrm{Q} M, \mrm{Q} I) \]
induced by $\mrm{Q}$. Finally given a morphism
$\phi: M \to I$ in $\cat{C}^{+}(\cat{Mod} A)$ we have
$\mrm{E}_{\msf{M}} \mrm{Q}(\phi) = \phi$.
\end{proof}

\section{Rigid Dualizing Complexes}
\label{sec3}

Dualizing complexes were introduced by Grothendieck \cite{RD}. 
The noncommutative variant was studied in \cite{Ye1}, and 
the notion of rigid dualizing complex is due to Van den Bergh 
\cite{VdB1}. Let us recall the definitions. From here to the end of 
the paper $\mbb{K}$ denotes a base field, and as before 
$\otimes = \otimes_{\mbb{K}}$.

An $A$-module $M$ is said to be {\em finite} if it is finitely 
generated. A homomorphism of rings $A \to B$ is called {\em finite} 
if $B$ is a finite $A$-module on both sides.
A $\mbb{K}$-algebra $A$ is called {\em affine} if it finitely 
generated.

\begin{dfn}[\cite{Ye1}, \cite{YZ2}] \label{dfn3.1}
Let $A$ be a left noetherian $\mbb{K}$-algebra and $B$ a right
noetherian $\mbb{K}$-algebra.
A complex $R \in \msf{D}^{\mrm{b}}(\cat{Mod} A \otimes B^{\mrm{op}})$
is called a {\em dualizing complex over $(A, B)$} if:
\begin{enumerate}
\rmitem{i} $R$ has finite injective dimension over $A$ and
$B^{\opn{op}}$.
\rmitem{ii} $R$ has finite cohomology modules over
$A$ and $B^{\opn{op}}$.
\rmitem{iii} The canonical morphisms
$B \to \opn{RHom}_{A}(R, R)$
in $\msf{D}(\cat{Mod} B^{\mrm{e}})$ and
$A \to \opn{RHom}_{B^{\opn{op}}}(R, R)$
in $\msf{D}(\cat{Mod} A^{\mrm{e}})$ are both isomorphisms.
\end{enumerate}
In case $A = B$, we shall say that $R$ is a dualizing complex
over $A$.
\end{dfn}

Condition (i) is equivalent to the existence of a quasi-isomorphism 
$R \to I$ in  
$\msf{C}^{\mrm{b}}(\cat{Mod} A \otimes B^{\mrm{op}})$
with each bimodule $I^{q}$ injective over $A$ and $B^{\opn{op}}$.

In this paper, whenever we mention a dualizing complex over
$(A, B)$ we implicitly assume that $A$ and $B^{\opn{op}}$ are
left noetherian $\mbb{K}$-algebras.

\begin{exa}
When $A$ is commutative and $R$ is a dualizing complex over $A$
consisting of central bimodules, then $R$ is a dualizing complex
in the sense of \cite[Section V.2]{RD}.
\end{exa}

\begin{dfn}[\cite{VdB1}]
Suppose $R$ is a dualizing complex over a noetherian $\mbb{K}$-algebra
$A$. If there is an isomorphism
\[ \phi: R \iso \opn{RHom}_{A^{\mrm{e}}}(A, R \otimes R) \]
in $\cat{D}(\cat{Mod} A^{\mrm{e}})$ we call the pair $(R, \phi)$
a {\em rigid dualizing complex}.
\end{dfn}

In the definition above $\opn{Hom}_{A^{\mrm{e}}}$ is with 
respect to the outside $A^{\mrm{e}}$-module structure of 
$R \otimes R$, and the isomorphism $\rho$ is with respect to the 
remaining inside $A^{\mrm{e}}$-module structure.

A rigid dualizing complex over $A$ is unique, up to an isomorphism
in $\cat{D}(\cat{Mod} A^{\mrm{e}})$, see
\cite[Proposition 8.2]{VdB1}. 

\begin{rem}
``Rigid dualizing complex'' is a relative notion, in the sense that 
it depends on the homomorphism $\mbb{K} \to A$. Cf.\ 
\cite[Example 3.13]{YZ2}.
\end{rem}

\begin{dfn}[\cite{YZ2}]
Suppose $A \to B$ is a finite homomorphism of $\mbb{K}$-algebras
and $(R_{A}, \phi_{A})$ and $(R_{B}, \phi_{B})$ are rigid
dualizing complexes over $A$ and $B$ respectively. A {\em rigid
trace} is a morphism $\opn{Tr}_{B / A} : R_{B} \to R_{A}$ in
$\cat{D}(\cat{Mod} A^{\mrm{e}})$ satisfying the two conditions
below.
\begin{enumerate}
\rmitem{i} $\opn{Tr}_{B / A}$ induces isomorphisms
\[ R_{B} \cong \opn{RHom}_{A}(B, R_{A}) \cong
\opn{RHom}_{A^{\mrm{op}}}(B, R_{A}) \]
in $\msf{D}(\cat{Mod} A^{\mrm{e}})$.
\rmitem{ii} The diagram
\[  \begin{CD}
R_{B} @> \phi_{B} >>
\opn{RHom}_{B^{\mrm{e}}}
(B, R_{B} \otimes R_{B}) \\
@V \opn{Tr} VV @V \opn{Tr} \otimes \opn{Tr} VV \\
R_{A} @> \phi_{A} >>
\opn{RHom}_{A^{\mrm{e}}}
(A, R_{A} \otimes R_{A})
\end{CD} \]
in $\msf{D}(\cat{Mod} A^{\mrm{e}})$ is commutative.
\end{enumerate}
\end{dfn}

According to \cite[Theorem 3.2]{YZ2} a rigid trace, if it exists,
is unique. Taking $A = B$ this implies that any two rigid dualizing
complexes $(R, \phi)$ and $(R', \phi')$ are {\em uniquely} isomorphic 
in $\msf{D}(\cat{Mod} A^{\mrm{e}})$, see \cite[Corollary 3.4]{YZ2}.
Often we shall omit explicit mention of the isomorphism $\phi$.

\begin{lem} \label{lem3.1}
Let $R$ be a dualizing complex over $(A, B)$, and let
$C :=$ \linebreak
$\opn{End}_{\msf{D}(\cat{Mod} A \otimes B^{\mrm{op}})}(R)$.
\begin{enumerate}
\item The left action of the center $\mrm{Z}(A)$ on $R$, 
and the right action of 
$\mrm{Z}(B)$ on $R$, induce isomorphisms of $\mbb{K}$-algebras
$\mrm{Z}(A) \cong C \cong \mrm{Z}(B)$.
These make $R$ into a complex of $C$-bimodules \tup{(}not necessarily 
central\tup{)}.
\item Let $M \in \msf{D}(\cat{Mod} A)$. Then the two $C$-module 
structures on $\opn{Ext}^{q}_{A}(M, R)$ coincide.
\item If $M \in \msf{D}(\cat{Mod} A \otimes B^{\mrm{op}})$ is 
$C$-central then the three $C$-module structures on 
$\opn{Ext}^{q}_{A \otimes B^{\mrm{op}}}(M, R \otimes R)$ coincide.
\item If $A = B$ and $R$ is rigid then the automorphism of 
$\mrm{Z}(A)$ in item \tup{(1)} is the identity.
\end{enumerate}
\end{lem}

\begin{proof}
The first item is a slight variation of \cite[Lemma 3.3]{YZ2} and 
\cite[Lemma 5.4]{Ye4}. In  item (2) the two actions of $C$ on 
$\opn{Ext}^{q}_{A}(M, R)$ correspond to the left action of $A$ on $R$ 
(and on $M$), and the right action of $B$ on $R$. Since 
$\opn{Ext}^{q}_{A}(M, R) = 
\opn{Hom}_{\msf{D}(\cat{Mod} A)}(M, R[q])$
these actions commute.  Likewise in item (3). 
Item (4) is \cite[Proposition 3.5]{YZ2}.
\end{proof}

\begin{lem} \label{lem3.3}
Let $A$ be a left noetherian $\mbb{K}$-algebra and 
$L \in \msf{D}^{-}_{\mrm{f}}(\cat{Mod} A)$.
\begin{enumerate}
\item Let $B$ be some $\mbb{K}$-algebra, let $N$ be a flat 
$B$-module and let
$M \in$ \linebreak 
$\msf{D}(\cat{Mod} (A \otimes B^{\mrm{op}}))$. Then
the canonical morphism
\[ \opn{RHom}_{A}(L, M) \otimes_{B} N \to 
\opn{RHom}_{A}(L, M \otimes_{B} N) \]
is an isomorphism.
\item Suppose $A \to A'$ a ring homomorphism such that 
$A'$ is a flat $A^{\mrm{op}}$-module. Let 
$M \in \msf{D}(\cat{Mod} A')$. Then the canonical morphism
\[ \opn{RHom}_{A}(L, M) \to 
\opn{RHom}_{A'}(A'\otimes_A L, M) \]
is an isomorphism. 
\end{enumerate}
\end{lem}

\begin{proof}
(1) Choose a quasi-isomorphism $P \to L$ where $P$ is a 
bounded above complex of finite free $A$-modules. Then the 
homomorphism of complexes
\[ \opn{Hom}_{A}(P, M) \otimes_{B} N \to 
\opn{Hom}_{A}(P, M \otimes_{B} N) \]
is bijective.

\medskip \noindent
(2) With $P \to L$ as above we get a free resolution
$A'\otimes_A P \to A'\otimes_A L$
as $A'$-modules, and 
\[ \opn{Hom}_{A}(P, M) \to \opn{Hom}_{A'}(A' \otimes_A P, M) \]
is bijective.
\end{proof}

The next theorem relates rigid dualizing complexes and central
localization. 

\begin{thm} \label{thm3.1}
Let $R$ be a dualizing complex over $(A, B)$, and identify  
$C \cong \mrm{Z}(A) \cong \mrm{Z}(B)$ as in Lemma 
\tup{\ref{lem3.1}}. Suppose $S \subset C$ is a multiplicatively
closed set, and let
$C_{S} := S^{-1} C$, $A_{S} := C_{S} \otimes_{C} A$ and
$B_{S} := C_{S} \otimes_{C} B$ be the localizations. Then:
\begin{enumerate}
\item The complex
\[ R_{S} := A_{S} \otimes_{A} R \otimes_{B} B_{S} \]
is a dualizing complex over $(A_{S}, B_{S})$.
\item If $A = B$, $R$ is a rigid dualizing complex over $A$, and 
$A^{\mrm{e}}$ is noetherian, then $R_{S}$ is a rigid dualizing 
complex over $A_S$.
\end{enumerate}
\end{thm}

\begin{proof}
(1) This is proved in a special case (when $A$ is commutative 
and $A_{S} = A_{\mfrak{p}}$ for a prime ideal $\mfrak{p}$)
in the course of the proof of \cite[Theorem 1.11(1)]{YZ2}; but 
the same proof works here too. Among other things one gets that 
$R_{S} \cong A_{S} \otimes_{A} R$ in 
$\msf{D}(\cat{Mod} A_{S} \otimes A^{\mrm{op}})$
and 
$R_{S} \cong R \otimes_{A} A_{S}$ in 
$\msf{D}(\cat{Mod} A \otimes A_{S}^{\mrm{op}})$

\medskip \noindent 
(2) We consider $R \otimes R$ as a left (resp.\ right) 
$A^{\mrm{e}}$-module via the outside (resp.\ inside) action. 
Since $A^{\mrm{e}}$ is noetherian and
$A^{\mrm{e}} \to (A_{S})^{\mrm{e}}$ is flat, by Lemma 
\ref{lem3.3}(1) we obtain an isomorphism
\[ \begin{aligned}
R_{S} 
& \cong 
R  \otimes_{A^{\mrm{e}}} (A_{S})^{\mrm{e}} \\
& \cong 
\opn{RHom}_{A^{\mrm{e}}}(A, R \otimes R) \otimes_{A^{\mrm{e}}}
(A_{S})^{\mrm{e}}  \\
& \cong 
\opn{RHom}_{A^{\mrm{e}}}(A, (R \otimes R) 
\otimes_{A^{\mrm{e}}} (A_{S})^{\mrm{e}})  
\end{aligned} \]
in $\msf{D}(\cat{Mod} (A_{S})^{\mrm{e}})$. 
Now
\[ (R \otimes R) \otimes_{A^{\mrm{e}}} (A_{S})^{\mrm{e}}
\cong R_S \otimes R_S \]
in $\msf{D}(\cat{Mod} (A^{\mrm{e}} \otimes (A_{S})^{\mrm{e}}))$. 
Finally using Lemma \ref{lem3.3}(2) we get
\[ \begin{aligned}
\opn{RHom}_{A^{\mrm{e}}}(A, R_S \otimes R_S) 
& \cong \opn{RHom}_{(A_S)^{\mrm{e}}}
((A_S)^{\mrm{e}} \otimes_A A, R_S \otimes R_S) \\
& \cong \opn{RHom}_{(A_S)^{\mrm{e}}}(A_S, R_S \otimes R_S) .
\end{aligned} \]
\end{proof}

If $M$ is a bimodule over a ring $A$ then the centralizer of $M$ 
is 
\[ \mrm{Z}_{A}(M) := \{ a \in A \mid a m = m a \text{ for all }
m \in M \} . \] 
Thus $\mrm{Z}_{A}(A) = \mrm{Z}(A)$. A ring homomorphism $A \to B$ 
is called {\em centralizing} if $B = A \cdot \mrm{Z}_{B}(A)$.
An {\em invertible bimodule} over $A$ is a bimodule $L$ such there 
exists another bimodule $L^{\vee}$ with 
$L \otimes_{A} L^{\vee} \cong L^{\vee} \otimes_{A} A \cong A$.
If $C$ is a commutative ring then a central invertible $C$-bimodule 
is the same as a projective $C$-module of rank $1$. 

\begin{prop} \label{prop3.3}
Suppose $C$ is a commutative affine $\mbb{K}$-algebra. Then $C$ 
has a rigid dualizing complex $R_{C}$ consisting of central 
bimodules. If $C$ is Cohen-Macaulay and equidimensional of dimension 
$n$ then we can choose
$R_{C} = \bsym{\omega}_{C}[n]$ where $\bsym{\omega}_{C}$ is a 
central bimodule, and if $C$ is Gorenstein then $\bsym{\omega}_{C}$ is 
invertible. 
\end{prop}

\begin{proof}
First assume 
$C = \mbb{K}[\bsym{t}] = \mbb{K}[t_{1}, \ldots, t_{n}]$,
a polynomial algebra. Then the bimodule $C$ is a dualizing complex. 
Because 
$\opn{Ext}^{n}_{C^{\mrm{e}}}(C, C^{\mrm{e}}) \cong C$
and 
$\opn{Ext}^{q}_{C^{\mrm{e}}}(C, C^{\mrm{e}}) = 0$
for $q \neq n$ it follows that the dualizing complex $C[n]$ is rigid.

Next take any affine algebra $C$. Choose a finite homomorphism 
$\mbb{K}[\bsym{t}] \to C$. Let $\mbb{K}[\bsym{t}] \to I$ be an 
injective resolution of the module $\mbb{K}[\bsym{t}]$ in 
$\msf{C}^{\mrm{b}}(\cat{Mod} C)$ and define
$R_{C} := \opn{Hom}_{\mbb{K}[\bsym{t}]}(C, I[n]) \in 
\msf{D}^{\mrm{b}}(\cat{Mod} C^{\mrm{e}})$.
So $R_{C}$ consists of central $C$-bimodules, and 
$R_{C} = \opn{RHom}_{\mbb{K}[\bsym{t}]}
(C, \mbb{K}[\bsym{t}][n])$.
According to the calculations in the proof of 
\cite[Proposition 5.7]{Ye4}, $R_{C}$ is a rigid dualizing complex.

Finally suppose $C$ is Cohen-Macaulay and equidimensional of dimension
$n$. Choose a noether normalization, that is a finite (and 
necessarily injective) homomorphism
$\mbb{K}[\bsym{t}] = \mbb{K}[t_{1}, \ldots, t_{n}] \to C$, 
cf.\ \cite[Theorem 13.3]{Ei}. According to \cite[Corollary 
18.17]{Ei},
$C$ is a projective $\mbb{K}[\bsym{t}]$-module. Hence
$\bsym{\omega}_{C} := \opn{Hom}_{\mbb{K}[\bsym{t}]}
(C, \mbb{K}[\bsym{t}])$
is a rigid dualizing complex. 
If moreover $C$ is Gorenstein then the bimodule $C$ is 
also a dualizing complex, and by the uniqueness of dualizing 
complexes over commutative algebras (cf.\ \cite[Theorem V.3.1]{RD})
we find that $\bsym{\omega}_{C}$ must be an invertible bimodule.  
\end{proof}

\begin{cor} \label{cor3.1}
Suppose $C$ is a commutative affine $\mbb{K}$-algebra, $R_{C}$ is 
a rigid dualizing complex over $C$ and $C \to A$ is a finite 
centralizing homomorphism of $\mbb{K}$-algebras. Then
$R_{A} := \opn{RHom}_{C}(A, R_{C})$
is a rigid dualizing complex over $A$.  
\end{cor}

\begin{proof}
Because of Proposition \ref{prop3.3} and the uniqueness of 
rigid dualizing complexes we may assume $R_{C}$ is a complex of 
central $C$-bimodules. Now proceed as in the proof of 
\cite[Proposition 5.7]{Ye4}.
\end{proof}

\begin{prop} \label{prop3.2}
Suppose $A$ is a noetherian affine $\mbb{K}$-algebra finite over 
its center and $A \to B$ is a finite centralizing homomorphism. 
Let $R_{A}$ and $R_{B}$ be rigid dualizing complexes over $A$ and 
$B$ respectively. The the rigid trace 
$\opn{Tr}_{B / A} : R_{B} \to R_{A}$
exists. 
\end{prop}

\begin{proof}
See \cite[Proposition 5.8]{Ye4}, noting that the morphism 
$\opn{Tr}_{B / A}$ constructed there satisfies axioms of rigid 
trace, as can be seen using the calculations done 
in the proof of \cite[Proposition 5.7]{Ye4}.
\end{proof}

\begin{rem}
An alternative approach to proving the last three results is via 
noetherian connected filtrations; see \cite[Theorem 7.16]{YZ2} and 
text prior to it.
\end{rem}

\begin{exa}
A rigid dualizing complex of a commutative $\mbb{K}$-algebra
$C$ need not be central. Let $R_{C}$ be as in Proposition 
\ref{prop3.3} above, and let $N^{0}$ be any non-central 
$C$-bimodule (e.g.\ $N^{0} = A^{\sigma}$, the twisted bimodule with 
$\sigma$ a nontrivial automorphism of $A$). Define the complex 
$N := (N^{0} \xrightarrow{=} N^{1})$. Then 
$R_{C} \cong R_{C} \oplus N$ in $\msf{D}(\cat{Mod} C^{\mrm{e}})$,
so the latter is a non-central rigid dualizing complex.
\end{exa}

\begin{exa} \label{exa3.1}
Assume $C$ is a smooth commutative $\mbb{K}$-algebra of relative
dimension $n$. Let $\Omega^{n}_{C / \mbb{K}}$ 
be the module of degree $n$ 
K\"{a}hler differentials. The canonical isomorphism (fundamental 
class of the diagonal)
\[ \Omega^{n}_{C / \mbb{K}} \cong \opn{Ext}^{n}_{C^{\mrm{e}}}
(C, \Omega^{2n}_{C^{\mrm{e}} / \mbb{K}}) \]
makes $\Omega^{n}_{C / \mbb{K}}[n]$ into a rigid dualizing 
complex. More generally for any $C$, if 
$\pi: \opn{Spec} C \to \opn{Spec} \mbb{K}$ is the structural 
morphism, then the twisted inverse image $\pi^{!} \mbb{K}$ of
\cite{RD} is the rigid dualizing complex of $C$.
\end{exa}

\section{Residual Complexes}
\label{sec4}

In this section we examine a refined notion of dualizing
complex, again generalizing from commutative algebraic geometry.
Some graded examples have been studied by Ajitabh \cite{Aj} and
the first author \cite{Ye2}. The main result here is Theorem
\ref{thm4.1}, which guarantees the existence of a residual complex.

Suppose $R$ is a dualizing complex over $(A, B)$ -- where $A$ and 
$B^{\mrm{op}}$ are left noetherian $\mbb{K}$-algebras --
and let $M$ be a finite $A$-module. The
grade of $M$ with respect to $R$ is
\[ j_{R; A}(M) := \inf \{ q \mid \opn{Ext}^{q}_{A}(M, R) \neq 0 \}
\in \mbb{Z} \cup \{\infty\} . \]
Similarly define $j_{R; B^{\opn{op}}}$ for a $B^{\opn{op}}$-module.

\begin{dfn}[\cite{Ye2}, \cite{YZ2}]  \label{dfn4.2}
Let $R$ be a dualizing complex over $(A, B)$.
We say that $R$ is an {\em Auslander dualizing complex} if it
satisfies the following two conditions:
\begin{enumerate}
\rmitem{i} For every finite $A$-module $M$,
integer $q$ and $B^{\opn{op}}$-submodule
$N \subset$ \linebreak
$\opn{Ext}^{q}_{A}(M, R)$,
one has $j_{R; B^{\opn{op}}}(N) \geq q$.
\rmitem{ii} The same holds after exchanging  $A$ and $B^{\opn{op}}$.
\end{enumerate}
\end{dfn}

\begin{dfn} \label{dfn4.1}
For a finite $A$-module $M$ the {\em canonical dimension} is
\[ \opn{Cdim}_{R; A} M := - j_{R; A}(M) . \]
\end{dfn}

It is known that if $R$ is an Auslander dualizing complex then
the canonical dimension $\opn{Cdim}_{R; A}$ is a spectral exact 
dimension function on $\cat{Mod} A$ (cf.\ Definition \ref{dfn1.1} 
and \cite[Theorem 2.10]{YZ2}). By symmetry there is a spectral exact
dimension function $\opn{Cdim}_{R; B^{\opn{op}}}$ on 
$\cat{Mod} B^{\opn{op}}$.

\begin{dfn} \label{dfn4.3}
A complex
$R \in \msf{C}^{\mrm{b}}(\cat{Mod} A \otimes B^{\mrm{op}})$
is called a {\em residual complex} over $(A, B)$ if the following
conditions are satisfied:
\begin{enumerate}
\rmitem{i} $R$ is a dualizing complex.
\rmitem{ii} Each bimodule $R^{-q}$ is an injective module over
$A$ and over $B^{\mrm{op}}$.
\rmitem{iii} $R$ is Auslander, and each bimodule $R^{-q}$ is
$\opn{Cdim}_{R; A}$-pure and $\opn{Cdim}_{R; B^{\mrm{op}}}$-pure
of dimension $q$ (Definition \ref{dfn2.3}).
\end{enumerate}
\end{dfn}

Let us denote by
$\mrm{Q}: \msf{C}(\cat{Mod} A \otimes B^{\mrm{op}}) \to
\msf{D}(\cat{Mod} A \otimes B^{\mrm{op}})$
the localization functor. If
$\msf{M} = \{ \msf{M}_{q}(\opn{Cdim}_{R; A}) \}$
then from Proposition \ref{prop2.7}(1) we see that
$\mrm{E}_{\msf{M}} \mrm{Q} R = R$ for a residual complex $R$.

A complex $I \in \cat{D}^{+}(\cat{Mod} A)$ is called a
minimal injective complex if the module $I^{q}$ is injective
and $\opn{Ker}(I^{q} \to I^{q + 1}) \subset I^{q}$
is essential, for all $q$.
Given $M \in \cat{D}^{+}(\cat{Mod} A)$, there is a quasi-isomorphism
$M \to I$ in  $\cat{C}^{+}(\cat{Mod} A)$, where $I$ is a
minimal injective complex.
Such $I$ is unique (up to a non-unique isomorphism),
and it is called the {\em minimal injective resolution} of $M$
(cf.\ \cite{Ye1} Lemma 4.2). 
If $M$ has finite injective dimension then $I$ is bounded.

\begin{dfn} \label{dfn4.4}
Let $R$ be an Auslander dualizing complex over $(A, B)$, and
let $I$ be the minimal injective resolution of $R$ in
$\cat{C}^{+}(\cat{Mod} A)$. Suppose each $A$-module $I^{-q}$
is $\opn{Cdim}_{R; A}$-pure of dimension $q$. Then we say $R$
has a {\em pure minimal injective resolution} over $A$. Likewise
for $B^{\mrm{op}}$.
\end{dfn}

According to \cite[Lemma 2.15]{Ye2}, a residual complex $R$ is a
minimal injective resolution of itself, on both sides. Thus $R$
has pure minimal injective resolutions.

\begin{prop} \label{prop4.1}
Suppose $R$ is an Auslander dualizing
complex over $(A, B)$ that has a pure minimal injective resolution
over $A$. Then the subcategories \linebreak
$\msf{M}_{q}(\opn{Cdim}_{R; A}) \subset \cat{Mod} A$
are stable for all $q$.
Likewise with $B^{\mrm{op}}$ and $A$ exchanged.
\end{prop}

\begin{proof}
We may use the proof of \cite[Proposition 2.7]{Ye2}.
\end{proof}

\begin{rem} \label{rem4.1}
If the subcategories 
$\msf{M}_{q}(\opn{Cdim}_{R; A}) \subset \cat{Mod} A$
are stable for all $q$ then $A$ is called a left pure algebra.
As shown in \cite{AjSZ}, many familiar algebras with Auslander
dualizing complexes do not admit
residual complexes -- indeed, are not even pure algebras
(cf.\ Example \ref{exa1.5}).
\end{rem}

\begin{lem} \label{lem4.2}
Let $R$ be a residual complex over $(A, B)$. Then the ring
homomorphisms \tup{(}left and right multiplication\tup{)}
\[ \mrm{Z}(A), \mrm{Z}(B^{\mrm{op}}) \to
\opn{Hom}_{\msf{C}(\cat{Mod} A \otimes B^{\mrm{op}})}(R, R) \]
are bijective.
\end{lem}

\begin{proof}
Since $R$ consists of injective $A$-modules,
$\opn{Hom}_{A}(R, R) = \opn{RHom}_{A}(R, R)$.
So
$\mrm{H}^{0} \opn{Hom}_{A}(R, R) = B^{\mrm{op}} \cdot 1_{R}
\cong B^{\mrm{op}}$.
By the purity assumption $\opn{Hom}_{A}(R, R)^{-1} = 0$, so
\[ \opn{Hom}_{\msf{C}(\cat{Mod} A)}(R, R) =
\mrm{H}^{0} \opn{Hom}_{A}(R, R) = B^{\mrm{op}} \cdot 1_{R}
\subset \opn{Hom}_{\msf{C}(\cat{Mod} \mbb{K})}(R, R) . \]
We see that
\[ \opn{Hom}_{\msf{C}(\cat{Mod} A \otimes B^{\mrm{op}})}(R, R) =
\mrm{Z}_{B^{\mrm{op}}}(B^{\mrm{op}} \cdot 1_{R}) =
\mrm{Z}(B^{\mrm{op}}) \cdot 1_{R} . \]
The equality with $\mrm{Z}(A) \cdot 1_{R}$ is proved symmetrically.
\end{proof}

\begin{thm} \label{thm4.1}
Suppose $R$ is an Auslander dualizing complex over $(A, B)$ that
has pure minimal injective resolutions over $A$ and
over $B^{\mrm{op}}$. Let
$\cat{M} = \{ \cat{M}_{q}(\opn{Cdim}_{R; A}) \}$ be the filtration
of $\cat{Mod} A$ determined by $R$ and let
$\mrm{E} := \mrm{E}_{\msf{M}}$ be the associated Cousin functor.
Then $\mrm{E} R$ is a residual complex, and there is a unique
isomorphism
\[ \phi: R \iso \mrm{Q} \mrm{E} R \text{ in }
\cat{D}(\cat{Mod} A \otimes B^{\mrm{op}})  \]
such that
\[ \mrm{E}(\phi): \mrm{E} R \to \mrm{E} \mrm{Q} \mrm{E} R
= \mrm{E} R \text{ in } \cat{C}(\cat{Mod} A \otimes B^{\mrm{op}}) \]
is the identity.
\end{thm}

\begin{proof}
If we decide to forget the $B^{\mrm{op}}$-module structure of $R$, we
may use the minimal injective resolution $I$ of $R$ as $A$-module to
compute
$\mrm{R} \Gamma_{\cat{M}^{p} /  \cat{M}^{p + 1}} R$.
By the purity assumption
\begin{equation} \label{eqn4.6}
\mrm{H}^{p + q} \mrm{R} \Gamma_{\cat{M}^{p} /  \cat{M}^{p + 1}} R
= \begin{cases}
I^{p} \quad & \text{if } q = 0 \\
0 & \text{otherwise} .
\end{cases}
\end{equation}
We conclude that $(\mrm{E}_{\cat{M}} R)^{-q} \cong I^{-q}$ as
$A$-modules, and so $(\mrm{E} R)^{-q}$ is an injective $A$-module,
$\opn{Cdim}_{R; A}$-pure of dimension $q$.

According to Proposition \ref{prop4.1} we see that the hypotheses of
Proposition \ref{prop2.2} hold with $M = R$,
$\cat{M} = \{ \cat{M}_{q}(\opn{Cdim}_{R; A}) \}$ and
$\cat{N} = \{ \cat{M}_{q}(\opn{Cdim}_{R; B^{\mrm{op}}}) \}$.
This tells us that
$\mrm{E}_{\cat{M}} R \cong \mrm{E}_{\cat{N}} R$
as complexes of bimodules. By the previous paragraph applied to
$B^{\mrm{op}}$ instead of $A$,
$(\mrm{E}_{\cat{N}} R)^{-q}$ is an injective $B^{\mrm{op}}$-module,
$\opn{Cdim}_{R; B^{\mrm{op}}}$-pure of dimension $q$.

Formula (\ref{eqn4.6}) says that Corollary \ref{cor2.2} holds here.
We deduce the existence of an isomorphism
$\phi': R \iso \mrm{Q} \mrm{E} R$ in
$\cat{D}(\cat{Mod} A \otimes B^{\mrm{op}})$. According to
Lemma \ref{lem4.2} there is some invertible element
$a \in \mrm{Z}(A)$ such that
$\mrm{E}(\phi'): \mrm{E} R \to \mrm{E} \mrm{Q} \mrm{E} R = \mrm{E} R$
is multiplication by $a$. The isomorphism
$\phi := a^{-1} \phi': R \to \mrm{Q} \mrm{E} R$ has the desired
property.
\end{proof}

Here is a class of algebras to which the results of this section 
apply. 
Recall that a ring $A$ is right bounded if every essential right 
ideal of $A$ contains an ideal which is essential as a right 
ideal. A ring $A$ is a right FBN (fully bounded noetherian) if 
$A$ is right noetherian and every prime factor ring of $A$ is 
right bounded. An {\em FBN ring} is a ring $A$ that is both
right and left FBN \cite[Chapter 8]{GW}.

A dualizing complex $R$ over two algebras $A$ and $B$ is called 
{\em weakly bifinite} if for every bimodule $M$ which 
is a subquotient of $A$, the bimodules 
$\opn{Ext}^{q}_{A}(M, R)$ are all finite on both sides; and the 
same is true with $A$ and $B^{\mrm{op}}$ interchanged.

An exact dimension function $\opn{dim}$, defined on 
$\cat{Mod} A$ and on $\cat{Mod} B^{\mrm{op}}$, is called 
{\em symmetric} if
$\opn{dim}_{A} M = \opn{dim}_{B^{\circ}} M$
for every bimodule $M$ finite on both sides.

\begin{thm} \label{thm4.2}
Let $A$ and $B$ be FBN $\mbb{K}$-algebras and let $R$ be an
Auslander dualizing complex over $(A,B)$ which
is weakly bifinite and such that $\opn{Cdim}_{R}$ is symmetric.
Then $R$ has pure minimal injective resolutions on both sides, and 
therefore the Cousin complex $\mrm{E} R$ \tup{(}notation as in 
Theorem \tup{\ref{thm4.1})} is a residual complex.
\end{thm}

For the proof we will need two lemmas and some notation.
Let $\mfrak{p}$ be a prime ideal of a left noetherian ring $A$. 
Write $S_{A / \mfrak{p}}(0)$ for the set of regular elements of 
$A / \mfrak{p}$. This is a denominator set in $A / \mfrak{p}$, and 
the ring of fractions 
$Q(\mfrak{p}) = \opn{Frac} A / \mfrak{p} = 
S_{A / \mfrak{p}}(0)^{-1} A / \mfrak{p}$ 
is simple artinian.

Given a finite $A$-module $M$, it's {\em reduced \tup{(}Goldie\tup{)} 
rank} at $\mfrak{p}$ is 
\[ \opn{rank}_{\mfrak{p}}(M) := 
\opn{length}_{Q(\mfrak{p})} Q(\mfrak{p}) \otimes_{A} M . \]
For $M = A / \mfrak{p}$ we write 
$\mrm{r}(\mfrak{p}) := \opn{rank}_{\mfrak{p}}(A / \mfrak{p})$. 
Let $J(\mfrak{p}) = J_{A}(\mfrak{p})$ be the indecomposable 
injective $A$-module with associated prime $\mfrak{p}$. 
The injective hull of $A / \mfrak{p}$ as $A$-module is then 
$J(\mfrak{p})^{\mrm{r}(\mfrak{p})}$.

Suppose $A$ is a prime ring. Recall that an element $m \in M$ is 
{\em torsion} if $a m = 0$ for some regular element $a \in A$. 
$M$ is a {\em torsion module} if all its elements are torsion; 
otherwise it is a {\em non-torsion module}. $M$ is {\em torsion-free}
if the only torsion element in it is $0$. 

\begin{lem} \label{lem4.4}
Suppose $A$ is a prime left noetherian ring, and $M, L$ are 
non-torsion $A$-modules, with $M$ finite. 
\begin{enumerate}
\item There is an injective homomorphism
$f : A \inj L^{r}$, where $r = \mrm{r}(0)$ is the Goldie rank.
\item Let $\opn{dim}$ be an exact dimension function on 
$\cat{Mod} A$. Then 
$\opn{dim} L = \opn{dim} A$.
\item There is a nonzero homomorphism 
$g : M \to L$. 
\end{enumerate}
\end{lem}

The proof of this lemma is standard, cf.\ 
\cite[Corollary 6.26(b)]{GW}.

The following is essentially proved in \cite[Lemma 2.3]{Br}. We 
state it for any minimal injective complex instead of a minimal 
resolution of a module.

\begin{lem} \label{lem4.3}
Suppose $A$ is a left noetherian ring. 
Let $I$ be a minimal injective complex of $A$-modules. Let $\mfrak{p}$
be a prime ideal of $A$, and let $\mu_{i}(\mfrak{p})$ be 
the multiplicity of $J(\mfrak{p})$ in $I^i$. Then
\begin{enumerate}
\item The image of the map 
$\partial^{i - 1}: \opn{Hom}_{A}(A / \mfrak{p}, I^{i-1}) \to 
\opn{Hom}_{A}(A / \mfrak{p}, I^{i})$
is a torsion $A/\mfrak{p}$-module.
\item $\mu_{i}(\mfrak{p}) = 
\opn{rank}_{\mfrak{p}}(\opn{Hom}_{A}(A / \mfrak{p}, I^{i})) =
\opn{rank}_{\mfrak{p}}(\opn{Ext}^{i}_{A}(A / \mfrak{p}, I))$.
\item{} Let $\mfrak{p}$ and $\mfrak{q}$ be two primes of $A$ and
$M$ an $A / \mfrak{p}$-$A / \mfrak{q}$-bimodule. Assume $M$ is 
nonzero, torsion-free as $(A / \mfrak{q})^{\mrm{op}}$-module, and 
finite non-torsion as $A / \mfrak{p}$-module. If $I^{i}$ contains 
a copy of $J(\mfrak{p})$, then $\opn{Ext}^{i}_{A}(M, I)$ is a 
non-torsion $A / \mfrak{q}$-module.
\end{enumerate}
\end{lem}

\begin{proof} (1) This is true because the kernel of this map is 
an essential submodule.

\medskip \noindent (2) The first equality is clear. The second follows 
from (1) because factoring out $A / \mfrak{p}$-torsion the complex 
$\opn{Hom}_{A}(A / \mfrak{p}, I)$ has zero coboundary maps. 

\medskip \noindent (3) 
By assumption there is a nonzero left ideal $L$ of $A / \mfrak{p}$ 
contained in $I^{i}$. Let 
$Z^{i} := \opn{Ker}(\partial^{i} : I^{i} \to I^{i+1})$.
Replacing $L$ by $L \cap Z^{i}$, we may assume 
$L \subset Z^{i}$. Since $L$ and $M$ are non-torsion 
$A / \mfrak{p}$-modules and $M$ is finite, Lemma \ref{lem4.4} says 
there is a nonzero map $f: M \to L$. 
We claim that $f$ is nonzero in 
$\opn{Ext}^{i}_{A}(M, I)$. Otherwise there is a map 
$g: M \to I^{i-1}$ such that $f = \partial^{i-1} g$. Let 
$M' := g(M) \subset I^{i-1}$. Then $M' \cap Z^{i-1}$ is essential
in $M'$ and hence $M' / (M'\cap Z^{i-1})$ is a torsion 
$A / \mfrak{p}$-module. Now $f(M)$ is a quotient of 
$M' / (M'\cap Z^{i-1})$, so it is a nonzero torsion 
$A / \mfrak{p}$-module. 
This contradicts the fact that any nonzero submodule of $L$ must be a 
torsion-free $A / \mfrak{p}$-module. Therefore we proved that $f$ 
is nonzero in $\opn{Ext}^{i}_{A}(M, I)$. 

Now let $a$ be a regular element of $A / \mfrak{q}$. Since 
$M \cong  M a$, $M / M a$ is a torsion $A / \mfrak{p}$-module. Hence 
$M a$ is not contained in the kernel of $f: M \to L$. This implies 
that $a f: M \to L$ is nonzero. By the claim we proved in the last 
paragraph, we see that $a f$ is nonzero in $\opn{Ext}^{i}_{A}(M, I)$.
So $f$ is non-torsion in $\opn{Ext}^{i}_{A}(M, I)$.
\end{proof}

\begin{proof}[Proof of Theorem \tup{\ref{thm4.1}}]
Let $I$ be the minimal injective resolution of $R$ as complex of 
$A$-modules. First we will show that each $I^{-i}$ is essentially 
pure of $\opn{Cdim}_{R; A} = i$, meaning that $I^{-i}$ contains 
an essential submodule that's pure of $\opn{Cdim}_{R; A} = i$.
It suffices to show that if $M$ is a $\opn{Cdim}_{R; A}$-critical 
submodule of $I^{-i}$, then $\opn{Cdim}_{R; A} M = i$. 

The critical module $M$ is uniform. Since $A$ is FBN the injective 
hull of $M$ is $J(\mfrak{q})$ for some prime ideal $\mfrak{q}$. 
Replacing $M$ by a nonzero submodule we can assume $M$ is a left 
ideal of $A / \mfrak{q}$, so it is a torsion-free $A / 
\mfrak{q}$-module. By Lemma 
\ref{lem4.3}(3), $E := \opn{Ext}^{-i}_{A}(A / \mfrak{q}, R)$ is a 
non-torsion $A / \mfrak{q}$-module. In particular $E \neq 0$
and hence $\opn{Cdim}_{R; A} A / \mfrak{q} \geq i$. By Lemma 
\ref{lem4.4} we get 
$\opn{Cdim}_{R; A} E = \opn{Cdim}_{R; A} A / \mfrak{q}$. 
 From the weakly bifinite hypothesis, $E$ is noetherian on both 
sides. Hence, by the symmetry of $\opn{Cdim}_{R}$, we have
$\opn{Cdim}_{R; A} E = \opn{Cdim}_{R; B^{\mrm{op}}} E$. 
According to \cite[Theorem 2.14]{YZ2} we have 
$\opn{Cdim}_{R; B^{\mrm{op}}} E \leq i$. We 
conclude that $\opn{Cdim}_{R; A} A / \mfrak{q} = i$. Again by 
Lemma \ref{lem4.4} we get 
$\opn{Cdim}_{R; A} M = \opn{Cdim}_{R; A} A / \mfrak{q} = i$.

Next we show that the $\opn{Cdim}_{R; A}$ is a constant on the cliques
of $A$. If there is a link $\mfrak{q} \rightsquigarrow \mfrak{p}$, 
then there is a nonzero
$A / \mfrak{q}$-$A / \mfrak{p}$-bimodule $M$ that is a subquotient 
of $A$ and is torsion free on both sides. 
By Lemma \ref{lem4.4} and the symmetry of $\opn{Cdim}_{R}$ we have
\begin{multline*}
\qquad \opn{Cdim}_{R; A} A / \mfrak{q} = \opn{Cdim}_{R; A} M
= \opn{Cdim}_{R; B^{\mrm{op}}} M \\
= \opn{Cdim}_{R; B^{\mrm{op}}} A / \mfrak{p} 
= \opn{Cdim}_{R; A} A / \mfrak{p} . \qquad
\end{multline*}

The FBN ring $A$ satisfies the second layer condition
\cite[4.3.14]{MR}. We know that $\opn{Cdim}_{R; A}$ is constant on cliques.
It follows from \cite[Proposition 4.3.13]{MR}
that an indecomposable injective $A$-module has pure 
$\opn{Cdim}_{R; A}$ (cf.\ \cite[Theorem 4.2]{AjSZ}). So the 
minimal injective resolution $I$ is pure. 

All the above works also for the minimal injective resolution of 
$R$ as a complex of $B^{\mrm{op}}$-modules.
\end{proof}

\begin{rem}
Let $A$ be a noetherian affine PI Hopf algebra over $\mbb{K}$ of 
finite injective dimension $n$. 
Brown and Goodearl \cite{BG} show that $A$ is Auslander-Gorenstein.
Using this one can show that the Auslander dualizing complex $A[n]$ 
is pre-balanced and has pure minimal injective resolutions (see 
\cite{YZ2}). According to Theorem \ref{thm4.1}, $\mrm{E} (A[n])$ 
is a residual complex. 
\end{rem}

\section{The Residue Complex of an Algebra}

In this section we define the residue complex of an algebra, 
combining the results if Sections \ref{sec3} and \ref{sec4}.
The main result here is Theorem \ref{thm5.1} which explains the 
functoriality of residue complexes. Here as before $\mbb{K}$ is 
the base field. 

If a $\mbb{K}$-algebra $A$ has a rigid dualizing complex $R$ 
(Definition \ref{dfn3.1}) that is Auslander (Definition 
\ref{dfn4.2}) then we shall usually write $\opn{Cdim}_{A}$ instead 
of $\opn{Cdim}_{R; A}$. This dimension function depends only on the  
$\mbb{K}$-algebra $A$.

\begin{dfn} \label{dfn5.1}
A {\em residue complex} over $A$ is a rigid dualizing complex
$(R, \phi)$ such that $R$ is also a residual complex (Definition
\ref{dfn4.3}).
\end{dfn}

The uniqueness of the residue complex will be made clear later in
this section (Corollary \ref{cor5.0}). In \cite{Ye2} the name
``strong residue complex'' was
used for the same notion (in the graded case).

Let $\cat{D}_{\mrm{f}}(\cat{Mod} A)$ denote the subcategory of
complexes with finite cohomology modules. The next result
will be used in the proof of Theorem \ref{thm5.1}.

\begin{prop}[Local Duality] \label{prop5.1}
Let $R$ be a dualizing complex over $(A,B)$ and
$\cat{M} \subset \cat{Mod} A$ a localizing subcategory.
Then there is a functorial isomorphism
\[ \mrm{R} \Gamma_{\cat{M}} M \cong
\opn{RHom}_{B^{\mrm{op}}}
(\opn{RHom}_{A}(M, R), \mrm{R} \Gamma_{\cat{M}} R) \]
for $M \in \cat{D}^{+}_{\mrm{f}}(\cat{Mod} A)$.
\end{prop}

\begin{proof}
Take a quasi-isomorphism $M \to I$ in $\cat{C}^{+}(\cat{Mod} A)$
with each $I^{q}$ injective over $A$, and a quasi-isomorphism
$R \to J$ in $\cat{C}^{\mrm{b}}(\cat{Mod} A \otimes B^{\mrm{op}})$,
where each $J^{q}$ is injective over $A$ and over $B^{\mrm{op}}$.
Write $D M := \opn{RHom}_{A}(M, R)$
and
$D^{\mrm{op}} N  := \opn{RHom}_{B^{\mrm{op}}}(N, R)$.
Using Lemma \ref{lem4.2} we get a commutative diagram in
$\cat{D}^{\mrm{b}}(\cat{Mod} A)$
\[ \choosegraphics{\UseTips \xymatrix{
{\mrm{R}} \Gamma_{\cat{M}} M \ar^{\alpha}[r] \ar^{=}[d] &
{\mrm{R}} \Gamma_{\cat{M}} D^{\mrm{op}} D M \ar^{\beta}[r]
\ar^{=}[d] &
\opn{RHom}_{B^{\mrm{op}}}(D M, \mrm{R} \Gamma_{\cat{M}} R)
\ar^{=}[d] \\
{\Gamma_{\cat{M}}} I \ar[r] &
{\Gamma_{\cat{M}}} \opn{Hom}_{B^{\mrm{op}}}(\opn{Hom}_{A}
(I, J), J)  \ar^{\gamma}[r] &
{\opn{Hom}}_{B^{\mrm{op}}}(\opn{Hom}_{A}
(I, J), \Gamma_{\cat{M}} J)
}} \]
with the bottom row consisting of morphisms in
$\cat{C}^{+}(\cat{Mod} A)$. The homomorphism $\gamma$ is actually
bijective. And since $M \in \cat{D}_{\mrm{f}}(\cat{Mod} A)$,
$M \to D^{\mrm{op}} D M$ is an isomorphism, and hence so is $\alpha$.
The isomorphism we want is $\beta \alpha$.
\end{proof}

\begin{rem}
The proposition above generalizes \cite[Theorem IV.6.2]{RD}
(resp.\ \cite[Theorem 4.18]{Ye1}), when $A$ is commutative local with
maximal ideal (resp.\ connected graded with augmentation ideal)
$\mfrak{m}$, and $\cat{M}$ is the category of $\mfrak{m}$-torsion
$A$-modules.
\end{rem}

Denote by
$\mrm{Q}: \cat{C}(\cat{Mod} A^{\mrm{e}}) \to
\cat{D}(\cat{Mod} A^{\mrm{e}})$
the localization functor.

\begin{thm} \label{thm5.1}
Let $A \to B$ be a finite centralizing homomorphism between
noetherian $\mbb{K}$-algebras. Suppose the two conditions below
hold.
\begin{enumerate}
\rmitem{i} There are  rigid dualizing complexes
$R_{A}$ and $R_{B}$ and the rigid trace morphism
$\opn{Tr}_{B / A} : R_{B} \to R_{A}$ exists.
\rmitem{ii} $R_{A}$ is an Auslander dualizing complex and it has
pure minimal injective resolutions on both sides.
\end{enumerate}
Then:
\begin{enumerate}
\item $R_{B}$ is an Auslander dualizing complex and it has
pure minimal injective resolutions on both sides.
\item Denote by
$\mrm{E}_{A}: \cat{D}^{+}(\cat{Mod} A^{\mrm{e}}) \to
\cat{C}^{+}(\cat{Mod} A^{\mrm{e}})$
and
$\mrm{E}_{B}: \cat{D}^{+}(\cat{Mod} B^{\mrm{e}}) \to
\cat{C}^{+}(\cat{Mod} B^{\mrm{e}})$
the Cousin functors associated to the dimension functions
$\opn{Cdim}_{R_{A}; A}$ and $\opn{Cdim}_{R_{B}; B}$ respectively.
Then
$\mrm{E}_{A} M \cong \mrm{E}_{B} M$
functorially for $M \in \cat{D}^{+}(\cat{Mod} B^{\mrm{e}})$.
\item Let $\mcal{K}_{A} := \mrm{E}_{A} R_{A}$ and
$\mcal{K}_{B} := \mrm{E}_{B} R_{B} \cong \mrm{E}_{A} R_{B}$
be the two residual complexes, so we have a morphism
$\mrm{E}_{A}(\opn{Tr}_{B / A}) : \mcal{K}_{B} \to \mcal{K}_{A}$
in $\cat{C}(\cat{Mod} A^{\mrm{e}})$.
Let
$\phi_{A}: R_{A} \iso \mrm{Q}  \mcal{K}_{A}$ and
$\phi_{B}: R_{B} \iso \mrm{Q}  \mcal{K}_{B}$ be
the isomorphisms from Theorem \tup{\ref{thm4.1}}. Then the diagram
\[ \choosegraphics{\UseTips \xymatrix{
R_{B} \ar^{\phi_{B}}[r] \ar^{\opn{Tr}}[d] &
\mrm{Q} \mcal{K}_{B} 
\ar^{\mrm{Q} \mrm{E}_{A}(\opn{Tr})}[d] \\
R_{A} \ar^{\phi_{A}}[r] &
\mrm{Q} \mcal{K}_{A} 
}} \]
in $\cat{D}(\cat{Mod} A^{\mrm{e}})$ is commutative.
\item $\mrm{E}_{A}(\opn{Tr}_{B / A})$ induces isomorphisms
\[ \mcal{K}_{B} \cong \opn{Hom}_{A}(B, \mcal{K}_{A}) \cong
\opn{Hom}_{A^{\mrm{op}}}(B, \mcal{K}_{A})  \]
in $\cat{C}(\cat{Mod} A^{\mrm{e}})$.
\end{enumerate}
\end{thm}

\begin{proof}
1. According to \cite[Proposition 3.9]{YZ2}
\[ R_{B} \cong \opn{RHom}_{A}(B, R_{A}) \cong
\opn{Hom}_{A}(B, \mcal{K}_{A})  \]
in $\cat{D}(\cat{Mod} (B \otimes A^{\mrm{op}}))$.
Therefore the complex $\opn{Hom}_{A}(B, \mcal{K}_{A})$ is an injective
resolution of $R_{B}$ in $\cat{K}^{+}(\cat{Mod} B)$.
Choose elements $b_{1}, \ldots, b_{r} \in \mrm{Z}_{B}(A)$ which
generate $B$ as an $A$-module. This gives rise to a surjection
$A^{r} \surj B$ of $A$-bimodules, and hence to an inclusion
$\opn{Hom}_{A}(B, \mcal{K}^{-q}_{A}) \subset (\mcal{K}_{A}^{-q})^{r}$.
So $\opn{Hom}_{A}(B, \mcal{K}_{A}^{-q})$ is
$\opn{Cdim}_{R_{A}; A}$-pure
of dimension $q$ as an $A$-module. By \cite[Proposition 3.9]{YZ2},
the dualizing complex $R_{B}$ is Auslander, and
$\opn{Cdim}_{R_{A}; A} M = \opn{Cdim}_{R_{B}; B} M$
for any $B$-module $M$. Thus $\opn{Hom}_{A}(B, \mcal{K}_{A}^{-q})$ is
$\opn{Cdim}_{R_{B}; B}$-pure as $B$-module.
We conclude that the injective resolution
$\opn{Hom}_{A}(B, \mcal{K}_{A})$ is
$\opn{Cdim}_{R_{B}; B}$-pure. But one easily sees that a pure
injective resolution must be minimal.
Symmetrically all the above applies to the right resolution
$\opn{Hom}_{A^{\mrm{op}}}(B, \mcal{K}_{A})$ of $R_{B}$.

\medskip \noindent 2.  Applying the functor $\mrm{E}_{B}$ to the
isomorphism $R_{B} \cong \opn{Hom}_{A}(B, \mcal{K}_{A})$ in
$\cat{D}(\cat{Mod} B \otimes A^{\mrm{op}})$,
and using the fact that $\opn{Hom}_{A}(B, \mcal{K}_{A})$
is a pure injective complex of $B$-modules, we obtain
$\mcal{K}_{B} \cong \opn{Hom}_{A}(B, \mcal{K}_{A})$ in
$\cat{C}(\cat{Mod} B \otimes A^{\mrm{op}})$.
Thus in particular $\opn{Hom}_{A}(B, \mcal{K}_{A})$ is a complex of
$B$-$B$-bimodules. By symmetry also
$\mcal{K}_{B} \cong \opn{Hom}_{A^{\mrm{op}}}(B, \mcal{K}_{A})$
as complexes of bimodules.

Denote by
$\msf{M}_{q}(A) := \{ M \in \cat{Mod} A \mid \opn{Cdim}_{R_{A}; A} M
\leq q \}$
and likewise $\msf{M}_{q}(B)$. We get filtrations
$\msf{M}(A) = \{ \msf{M}_{q}(A) \}$ and
$\msf{M}(B) = \{ \msf{M}_{q}(B) \}$. Since
$\Gamma_{\msf{M}_{q}(A)} \mcal{K}_{A} = \mcal{K}_{A}^{\geq -q}$,
Proposition \ref{prop5.1} tells us that
\[ \begin{split}
\mrm{R} \Gamma_{\cat{M}_{q}(A)} M & \cong
\opn{Hom}_{A^{\mrm{op}}}(\opn{Hom}_{A}(M, \mcal{K}_{A}),
\mcal{K}_{A}^{\geq -q}) \\
& \cong
\opn{Hom}_{{B}^{\mrm{op}}}(\opn{Hom}_{B}(M, \mcal{K}_{B}),
\mcal{K}_{B}^{\geq -q})
\cong
\mrm{R} \Gamma_{\cat{M}_{q}(B)} M
\end{split} \]
functorially for $M \in \cat{D}^{+}_{\mrm{f}}(\cat{Mod} B)$.
In particular
$\mrm{H}^{p}_{\cat{M}_{q}(B)} M \cong
\mrm{H}^{p}_{\cat{M}_{q}(A)} M$
functorially for finite $B$-modules $M$. Passing to
direct limits (using Proposition \ref{prop1.6}) 
this becomes true for all $B$-modules. Hence if $M$
is an $\msf{M}(B)$-flasque $B$-module, it is also
$\msf{M}(A)$-flasque. By Proposition \ref{prop2.6} it follows that
$\mrm{E}_{A} M \cong \mrm{E}_{B} M$
functorially for $M \in \cat{D}^{+}(\cat{Mod} B^{\mrm{e}})$.

\medskip \noindent 3.
Next we analyze the morphism
$\mrm{Q} \mrm{E}_{A}(\opn{Tr}_{B / A}) \in
\opn{Hom}_{\cat{D}(\cat{Mod} A^{\mrm{e}})}(R_{B}, R_{A})$.
By \cite[Lemma 3.3]{YZ2},
\[ \opn{Hom}_{\cat{D}(\cat{Mod} A^{\mrm{e}})}(R_{B}, R_{A}) =
\mrm{Z}_{B}(A) \cdot \opn{Tr}_{B / A} , \]
so
$\mrm{Q} \mrm{E}_{A}(\opn{Tr}_{B / A}) =
b \cdot \opn{Tr}_{B / A}$
for some (unique) $b \in \mrm{Z}_{B}(A)$.
We shall prove that $b = 1$.
If we forget the $A^{\mrm{op}}$-module structure, then
\[ \opn{Hom}_{\cat{D}(\cat{Mod} A)}(\mrm{Q} \mcal{K}_{B},
\mrm{Q} \mcal{K}_{A}) =
B \cdot \opn{Tr}_{B / A}. \]
 From Proposition \ref{prop2.7} we get that there is a bijection
\[ \opn{Hom}_{\cat{D}(\cat{Mod} A)}(\mrm{Q} \mcal{K}_{B},
\mrm{Q} \mcal{K}_{A}) \cong
\opn{Hom}_{\cat{C}(\cat{Mod} A)}(\mcal{K}_{B}, \mcal{K}_{A}) \]
induced by $\mrm{E}$, and the inverse is induced by $\mrm{Q}$.
Hence we obtain
\[ \opn{Tr}_{B / A} = \mrm{Q} \mrm{E}_{A}(\opn{Tr}_{B / A})
\in \opn{Hom}_{\cat{D}(\cat{Mod} A)}(\mrm{Q} \mcal{K}_{B},
\mrm{Q} \mcal{K}_{A}) . \]
This implies that $b = 1$.

\medskip \noindent 4.\
By part 3, if we apply the functor $\opn{Hom}_{A}(B, -)$ to the
homomorphism of complexes
$\mrm{E}_{A}(\opn{Tr}_{B / A}) : \mcal{K}_{B} \to \mcal{K}_{A}$
we get a quasi-isomorphism
$\mcal{K}_{B} \to \opn{Hom}_{A}(B, \mcal{K}_{A})$.
But these are minimal injective complexes of $B$-modules, so it must
actually be an isomorphism of complexes. By symmetry also
$\mcal{K}_{B} \to \opn{Hom}_{A^{\mrm{op}}}(B, \mcal{K}_{A})$
is an isomorphism of complexes.
\end{proof}

\begin{cor}[Uniqueness of Residue Complex] \label{cor5.0}
Suppose $(\mcal{K}_{A}, \phi)$ and $(\mcal{K}'_{A}, \phi')$ are
two residue complexes over $A$. Then there is a unique isomorphism
$\tau: \mcal{K}'_{A} \iso \mcal{K}_{A}$
in $\msf{C}(\cat{Mod} A^{\mrm{e}})$
that's compatible with $\phi'$ and $\phi$, i.e.\ a rigid trace.
\end{cor}

\begin{proof}
Write $B := A$ and
$(\mcal{K}_{B}, \phi_{B}) := (\mcal{K}'_{A}, \phi')$.
By \cite[Theorem 3.2]{YZ2} we get a unique isomorphism
$\opn{Tr}_{B / A}: \mcal{K}_{B} \iso \mcal{K}_{A}$
in $\msf{D}(\cat{Mod} A^{\mrm{e}})$ that's a rigid trace.
According to part 3 of the theorem above,
$\tau := \mrm{E}_{A}(\opn{Tr}_{B / A})$
satisfies $\mrm{Q}(\tau) = \opn{Tr}_{B / A}$, so it too is a
rigid trace.
\end{proof}

\begin{cor} \label{cor5.1}
If in the previous theorem $B = A / \mfrak{a}$ for some 
ideal $\mfrak{a}$ then there is equality
\[ \opn{Hom}_{A}(A / \mfrak{a}, \mcal{K}_{A}) =
\opn{Hom}_{A^{\mrm{op}}}(A / \mfrak{a}, \mcal{K}_{A}) \subset
\mcal{K}_{A} . \]
\end{cor}

\begin{proof}
By part 4 of the theorem we get an isomorphism of
$B \otimes A^{\mrm{op}}$-modules
$\mcal{K}_{B}^{q} \cong \opn{Hom}_{A}(B, \mcal{K}_{A}^{q})
\subset \mcal{K}_{A}^{q}$
for every $q$. This implies that
$\opn{Hom}_{A}(B, \mcal{K}_{A}^{q})$ is annihilated by $\mfrak{a}$
on the right too, and hence
$\opn{Hom}_{A}(B, \mcal{K}_{A}) \subset
\opn{Hom}_{A^{\mrm{op}}}(B, \mcal{K}_{A})$.
By symmetry there is equality.
\end{proof}

\begin{rem}
Corollary \ref{cor5.1} is pretty surprising. The ideal $\mfrak{a}$
will in general not be generated by central elements. On the other
hand the centralizer $\mrm{Z}_{A}(\mcal{K}_{A}) = \mrm{Z}(A)$.
So there is no obvious reason for the left annihilator of
$\mfrak{a}$ in $\mcal{K}_{A}$ to coincide with the right
annihilator.
\end{rem}

\begin{cor} \label{cor5.2}
Let $A \to B$ and $B \to C$ be finite centralizing homomorphisms.
Assume the hypotheses of Theorem \tup{\ref{thm5.1}}, and also that 
the rigid dualizing complex $R_{C}$ and the rigid trace
$\opn{Tr}_{C / B}$ exist. Then
\[ \mrm{E}_{A}(\opn{Tr}_{C / A}) =
\mrm{E}_{A}(\opn{Tr}_{B / A}) \, \mrm{E}_{B}(\opn{Tr}_{C / B}) :
\mcal{K}_{C} \to \mcal{K}_{A} . \]
\end{cor}

\begin{proof}
By \cite[Corollary 3.8]{YZ2} the morphism
$\opn{Tr}_{C / A} := \opn{Tr}_{B / A} \opn{Tr}_{C / B}$
is a rigid trace. According to Theorem \ref{thm5.1} the residue 
complex $\mcal{K}_{C} = \mrm{E}_{B} R_{C}$ exists, and 
$\mrm{E}_{B}(\opn{Tr}_{C / B}) = \mrm{E}_{A}(\opn{Tr}_{C / B})$.
\end{proof}

Here are a few examples of algebras with residue complexes. 

\begin{exa} \label{exa5.1}
If $A$ is a commutative affine (i.e.\ finitely generated) 
$\mbb{K}$-algebra and $R_{A}$ is 
its rigid dualizing complex then the complex 
$\mcal{K}_{A} := \mrm{E} R_{A}$ is a 
residue complex. It consists of central bimodules, and is
the residue complex of $A$ also in the sense of \cite{RD}; 
cf.\ Example \ref{exa3.1}. For a finite homomorphism $A \to B$ 
of commutative algebras the trace morphism $\opn{Tr}_{B / A}$
coincides with that of \cite{RD}. 
\end{exa}

\begin{exa} \label{exa5.2}
Consider a commutative artinian local $\mbb{K}$-algebra $A$ whose 
residue field $A / \mfrak{m}$ is finitely generated over $\mbb{K}$
(i.e.\ $A$ is residually finitely generated). Then
$A \cong \opn{Frac} A_{0}$, the ring of fractions of some commutative
affine $\mbb{K}$-algebra $A_{0} \subset A$, and by Theorem 
\ref{thm3.1} $A$ has a rigid dualizing complex 
$R_{A} \cong A \otimes_{A_{0}} R_{A_{0}}$.
Because of the uniqueness of dualizing complexes for commutative 
algebras, $R_{A} \cong \mcal{K}(A)[n]$ where 
$\mcal{K}(A) := \mrm{H}^{-n}_{\mfrak{m}} R_{A}$
is an injective hull of $A / \mfrak{m}$ and 
$n = \opn{dim} A_{0} = \opn{tr{.}deg}_{\mbb{K}}(A / \mfrak{m})$.

If $A \to B$ is a finite homomorphism of such artinian algebras 
then the rigid trace
$\opn{Tr}_{B / A} : \mcal{K}(B)[n] \to \mcal{K}(A)[n]$
exists. 

Now if $A$ is a residually finitely generated commutative noetherian
complete local $\mbb{K}$-algebra we can define 
$\mcal{K}(A) := \varinjlim \mcal{K}(A / \mfrak{m}^{i})$. 
The functorial $A$-module $\mcal{K}(A)$ is called the {\em dual 
module} of $A$. Cf.\ \cite{Ye3} and \cite{Ye5} for alternative 
approaches, applications and references to other related work.
\end{exa}

\begin{exa}
If $A$ is a noetherian affine $\mbb{K}$-algebra finite over its 
center $C$ and $\mcal{K}_{C}$ is the residue complex of $C$ then 
$\mcal{K}_{A} := \opn{Hom}_{C}(A, \mcal{K}_{C})$
is the residue complex of $A$. If $A \to B$ is a finite centralizing 
homomorphism then the rigid trace
$\opn{Tr}_{B / A} : \mcal{K}_{B} \to \mcal{K}_{A}$ 
is gotten by applying $\opn{Hom}_{C}(-, \mcal{K}_{C})$
to $A \to B$. We see that the theory of residual complexes for 
algebras finite over the center is very close to the commutative 
theory.
\end{exa}

\begin{prop} \label{prop5.2}
The first Weyl algebra over the field $\mbb{C}$ has a residue 
complex.
\end{prop}

\begin{proof}
Recall that the first Weyl algebra is 
$A := \mbb{C} \bra{x, y} / (y x - x y - 1)$.
According to \cite[Example 6.20]{YZ2}, $R_{A} := A[2]$ is a rigid 
Auslander dualizing complex over $A$, and $\opn{Cdim} = \opn{GKdim}$. 
The ring of fractions $Q = \opn{Frac} A$ is a division ring, and the
global dimension of $A$ is $1$. Therefore the minimal injective
resolution of $A$ in $\cat{Mod} A$ is
$0 \to A \to I^{0} \to I^{1} \to 0$ and $I^{0} \cong Q$.
We see that $I^{0}$ is pure of 
$\opn{GKdim} I^{0} = \opn{GKdim} A = 2$.
Since $I^{1}$ is a torsion $A$-module we get
$\opn{GKdim} I^{1} \leq 1$; but since there are no $A$-modules of 
$\opn{GKdim} = 0$ it follows that $I^{1}$ is pure of
$\opn{GKdim} = 1$.
So $A$ has a pure injective resolution on the left. 
The same is true on the right too. 
We see that $R_{A}$ has pure minimal injective resolutions on
both sides, so $\mcal{K}_{A} := \mrm{E} R_{A}$ is a residue complex.
Moreover $\mcal{K}_{A}^{-2} = Q$ and 
$\mcal{K}_{A}^{-1} = Q / A$.
\end{proof}

\begin{prop} \label{prop5.3}
Let $\mfrak{g}$ be a nilpotent $3$-dimensional Lie 
algebra over $\mbb{C}$ and $A := \mrm{U}(\mfrak{g})$ the universal 
enveloping algebra. Then $A$ has a residue complex.
\end{prop}

\begin{proof}
We may assume $A$ is not commutative, so $A$ is generated by 
$x, y, z$ with $z$ central and $[x, y] = z$. By 
\cite[Proposition 6.18]{YZ2} and \cite[Theorem A]{Ye5}
the complex $R_{A} := A[3]$ is a rigid Auslander dualizing complex, 
and $\opn{Cdim} = \opn{GKdim}$. Consider a minimal injective 
resolution 
$0 \to A \to I^{0} \to I^{1} \to I^{2} \to I^{3} \to 0$
of $A$ in  $\cat{Mod} A$. For any $\lambda \in \mbb{C}$ consider 
the ideal $\mfrak{a} = A \cdot (z - \lambda)$. The 
localizing subcategory 
$\msf{M}_{(z - \lambda)} = \msf{M}_{\mfrak{a}} \subset \cat{Mod} A$ 
is stable (cf.\ Example \ref{exa1.4}). We get a direct sum 
decomposition of $A$-modules indexed by $\opn{Spec} \mbb{C}[z]$:
\[ I^{q} \cong \Bigl( \bigoplus_{\lambda \in \mbb{C}} 
\Gamma_{(z - \lambda)} I^{q} \Bigr) \oplus 
\bigl( \mbb{C}(z) \otimes_{\mbb{C}[z]} I^{q} \bigr) . \]

For any $\lambda \in \mbb{C}$,  
$\opn{RHom}_{A}(A / (z - \lambda), R_{A})$
is the rigid dualizing complex of \linebreak 
$A / (z - \lambda)$, and
$\opn{Hom}_{A}(A / (z - \lambda), I[3])$
is its minimal injective resolution as complex of left modules. 
Since the algebra $A / (z - \lambda)$ is isomorphic to either the 
commutative polynomial algebra ($\lambda = 0$) or to the first 
Weyl algebra, we see that 
$\opn{Hom}_{A}(A / (z - \lambda), I^{q})$
is pure of $\opn{GKdim} = 3 - q$, for $1 \leq q \leq 3$.

Fix $q$ and $\lambda$. Introduce a filtration $F$ on 
$N := \Gamma_{(z - \lambda)} I^{q}$
by $F^{-j} N :=$ \linebreak
$\opn{Hom}_{A}(A / (z - \lambda)^{j}, I^{q})$.
Then for $j \geq 1$ multiplication by $z - \lambda$ is a bijection
$\opn{gr}_{F}^{-j - 1} N \iso \opn{gr}_{F}^{-j} N$. 
It follows that $\opn{gr}_{F}^{-j} N$ is pure of 
$\opn{GKdim} = 3 - q$. Therefore also 
$N = \bigcup F^{-j} N$ is pure of $\opn{GKdim} = 3 - q$.

The direct sum complement $\mbb{C}(z) \otimes_{\mbb{C}[z]} I^{q}$
is a $B$-module, where $B := \mbb{C}(z) \otimes_{\mbb{C}[z]} A$. 
In fact $\mbb{C}(z) \otimes_{\mbb{C}[z]} I$ is a minimal injective 
resolution of $B$ in $\cat{Mod} B$.
But $B$ is isomorphic the first Weyl algebra over the field 
$\mbb{C}(z)$. Therefore $I^{0} \cong \opn{Frac} B$ is pure of 
$\opn{GKdim} = 3$, and 
$\mbb{C}(z) \otimes_{\mbb{C}[z]} I^{1}$ 
is pure of $\opn{GKdim} = 2$.

We conclude that each $I^{q}$ is pure of 
$\opn{GKdim} = 3 - q$. By symmetry the same is true on the right 
too. So $\mcal{K}_{A} := \mrm{E} R_{A}$ is a residue complex over $A$.
\end{proof}

\begin{rem}
There are nilpotent Lie algebras $\mfrak{g}$ such that 
$\mrm{U}(\mfrak{g})$ does not have a residue complex. Indeed one 
can find such a Lie algebra with a surjection from 
$A = \mrm{U}(\mfrak{g})$ to the second Weyl algebra $B$. 
$B$ is not pure, so it does not have a residue complex. Hence by 
Theorem \ref{thm5.1}, $A$ does not have a residue complex.
\end{rem}

\begin{exa} \label{exa5.4}
Let $A$ be a $3$-dimensional Sklyanin algebra over the 
algebraically closed field $\mbb{K}$.
The whole apparatus of Cousin functors can be implemented also in the 
$\mbb{Z}$-graded module category $\cat{GrMod} A$-- actually this was
already done in \cite{Ye2} -- and in particular 
Theorem \ref{thm4.1} is true in the graded sense. According to 
\cite{Aj} the minimal injective resolutions of $A$ in 
$\cat{GrMod} A$ and in $\cat{GrMod} A^{\mrm{op}}$ are  
pure. On the hand the balanced dualizing complex, which is also 
rigid in the graded sense, is $R_{A} = \bsym{\omega}_{A}[3]$ where
$\bsym{\omega}_{A} = A^{\sigma}$ for some automorphism $\sigma$. 
We conclude that $\mcal{K}_{A} := \mrm{E} R_{A}$ is a graded residue 
complex over $A$. Note that this result was proved in \cite{Ye2}
by a direct (and rather involved) calculation of \"{O}re localizations 
with respect to $\sigma$-orbits in the elliptic curve associated 
to $A$.
\end{exa}

\begin{que}
In case the rigid Auslander dualizing complex $R$ exists but there 
is no residue complex (e.g.\ $A = \mrm{U}(\mfrak{sl}_{2})$),
is it still true that 
$\mrm{E}_{A} R \cong \mrm{E}_{A^{\mrm{op}}} R$?
What can be said about this complex? 
\end{que}

In the following section we will discuss residue complexes over PI 
algebras in detail.

\section{The Residue Complex of a PI Algebra}

In this section we look at an affine noetherian PI algebra $A$ 
over the base field $\mbb{K}$. We show that -- under a certain 
technical assumption -- such an algebra $A$ has a residue complex 
$\mcal{K}_{A}$. Furthermore in Theorem \ref{thm6.2} we give a 
detailed description of  the structure of $\mcal{K}_{A}$. 
The material on PI rings needed here can be found 
in \cite[Section 13]{MR}. 

\begin{prop} \label{prop6.2}
Suppose $A$ is an affine prime PI $\mbb{K}$-algebra with center 
$C$. Then there is a nonzero element $s \in C$
such that the localization $A_{s}$ is an Azumaya algebra over its 
center $\mrm{Z}(A_{s}) = C_{s}$, and $C_{s}$ is a regular 
commutative affine $\mbb{K}$-algebra.
\end{prop}

\begin{proof}
By the Artin-Procesi Theorem \cite[Theorem 13.7.14]{MR} 
and \cite[Proposition 13.7.4]{MR} we may 
find $s_{1} \in C$, $s_{1} \neq 0$ such that 
$A_{s_{1}}$ is an Azumaya algebra over its center 
$C_{1} := \mrm{Z}(A_{s_{1}})$. 
The commutative prime $\mbb{K}$-algebra $C_{1}$ is affine, and 
hence by \cite[page 246, Theorem 73]{Mat} there is a nonzero
element $s_{2} \in C_{1}$ such that the localization 
$C_{2} := (C_{1})_{s_{2}}$ is regular. By Posner's Theorem 
\cite[Theorem 13.6.5]{MR} the fraction fields coincide:
$\opn{Frac} C = \opn{Frac} C_{2}$. Because $C_{2}$ is affine
we may find $s \in C$ (the product of the denominators of a 
finite set of $\mbb{K}$-algebra generators of $C_{2}$)
such that $C_{s} = (C_{2})_{s}$. Hence $C_{s}$ is also regular, 
affine over $\mbb{K}$, and $A_{s}$ is Azumaya with center $C_{s}$. 
\end{proof}

\begin{thm} \label{thm6.3}
Let $A$ be an affine prime noetherian PI $\mbb{K}$-algebra with 
center $C$ and Gelfand-Kirillov dimension $\opn{GKdim} A = n$.
Assume $A$ has a rigid dualizing complex $R_{A}$.
\begin{enumerate}
\item Let $s \in C$ be a nonzero element such that the localization 
$A_{s}$ is an Azumaya 
algebra with center $C_{s}$, and $C_{s}$ is a regular affine 
$\mbb{K}$-algebra. Then there is an isomorphism
\[ C_{s} \otimes_{C} R_{A} \cong \bsym{\omega}_{C_{s}}[n] 
\otimes_{C} A \]
in  $\msf{D}(\cat{Mod} A^{\mrm{e}})$, where $\bsym{\omega}_{C_{s}}$ 
is a projective $C_{s}$-module of rank $1$.
\item Let $K := \opn{Frac} C$ and $Q := \opn{Frac} A$. Then
\[ K \otimes_{C} R_{A} \cong Q[n] \]
in  $\msf{D}(\cat{Mod} A^{\mrm{e}})$.
\end{enumerate}
\end{thm}

\begin{proof}
(1) By \cite[Proposition 4.4]{ASZ} the algebra 
$A^{\mrm{e}}$ is noetherian. 
Therefore according to Theorem \ref{thm3.1} the complex 
$R_{A_{s}} := A_{s} \otimes_{A} R_{A} \otimes_{A} A_{s}$
is a rigid dualizing complex over $A_{s}$. Moreover, 
\[ R_{A_{s}} \cong R_{A} \otimes_{C} C_{s} 
\cong C_{s} \otimes_{C} R_{A} \]
in $\msf{D}(\cat{Mod} A^{\mrm{e}})$.

By \cite[Proposition 8.2.13]{MR} we have
$\opn{GKdim} A = \opn{GKdim} A_{s}$, and hence
$n = \opn{GKdim} C_{s} = \opn{Kdim} C_{s}$ (Krull dimension). 
According to Proposition \ref{prop3.3} the rigid dualizing 
complex of $C_{s}$ is $\bsym{\omega}_{C_{s}}[n]$ with 
$\bsym{\omega}_{C_{s}}$ a projective $C_{s}$-module of rank $1$. 
 From Corollary \ref{cor3.1} we see that 
$\opn{Hom}_{C_{s}}(A_{s}, \bsym{\omega}_{C_{s}})$
is a rigid dualizing complex over $A_{s}$. Finally the reduced 
trace $A_{s} \to C_{s}$ induces a bimodule isomorphism
$A_{s} \cong \opn{Hom}_{C_{s}}(A_{s}, C_{s})$.
Therefore 
\[ R_{A_{s}} \cong \bsym{\omega}_{C_{s}} \otimes_{C_{s}} A_{s}
\cong \bsym{\omega}_{C_{s}} \otimes_{C} A \]
in $\msf{D}(\cat{Mod} A^{\mrm{e}})$.

\medskip \noindent 
(2) Follows from (1).
\end{proof}

Recall that a {\em connected graded} $\mbb{K}$-algebra is an
$\mbb{N}$-graded algebra
$A = \bigoplus_{i \in \mbb{N}} A_{i}$
such that $A_{0} = \mbb{K}$ and
$\opn{rank}_{\mbb{K}} A_{i} < \infty$
for all $i$. By a filtration of $A$ we mean an ascending
filtration $F = \{ F_{i} A \}_{i \in \mbb{Z}}$ by
$\mbb{K}$-modules such that
$F_{i} A \cdot F_{j} \subset F_{i + j} A$.
The associated graded algebra is denoted by $\opn{gr}^{F} A$.

\begin{dfn} \label{dfn6.1}
A {\em noetherian connected filtration} of a $\mbb{K}$-algebra $A$
is a filtration $F$ such that $\opn{gr}^{F} A$ is a noetherian
connected graded $\mbb{K}$-algebra.
\end{dfn}

In \cite[Definition 6.1]{YZ2} the condition was that the Rees algebra
$\opn{Rees}^{F} A$ should be a noetherian connected graded
$\mbb{K}$-algebra; but as mentioned there the two conditions
are in fact equivalent.

It is not hard to see that if $A$ admits a noetherian connected
filtration then $A$ itself is noetherian and affine over
$\mbb{K}$.

\begin{rem}
If $A$ is affine and finite over its center then it admits a
noetherian connected filtration (see \cite[Example 6.14]{YZ2});
but this case is in a sense too easy. There are known examples of
PI algebras not finite over their centers that admit noetherian
connected filtrations (e.g.\ \cite[Example 6.15]{YZ2}), 
and for a long time it was an open problem
whether they all do. The first counterexample was recently discovered
by Stafford \cite{SZ}.
\end{rem}

The notions of symmetric dimension function and weakly bifinite 
dualizing complex were defined just before Theorem
\ref{thm4.2}.

\begin{prop} \label{prop6.3}
Suppose $A$ is a PI algebra admitting a noetherian connected 
filtration. Then $A$ has an Auslander rigid dualizing complex $R$,
the canonical dimension $\opn{Cdim} = \opn{Cdim}_{R}$ 
is symmetric, and $R$ is weakly bifinite. 
\end{prop}

\begin{proof}
Let $F$ be a noetherian connected filtration of $A$. Then 
$\opn{gr}^{F} A$ is a noetherian connected graded PI 
$\mbb{K}$-algebra. By \cite[Corollary 6.9]{YZ2} $A$ has an 
Auslander dualizing complex $R$, and 
$\opn{Cdim}_{R} = \opn{GKdim}$ (Gelfand-Kirillov dimension) on the 
categories $\cat{Mod} A$ and $\cat{Mod} A^{\mrm{op}}$. Since 
$\opn{GKdim}$ is symmetric (see 
\cite[Proposition 8.3.14(ii)]{MR}),
so is $\opn{Cdim}_{R}$.

Now take a bimodule $M$ that's a subquotient of $A$. Then $M$ 
admits a two-sided good filtration $F$ (i.e.\ $\opn{gr}^{F} M$ is 
a finite module over $\opn{gr}^{F} A$ on both sides), and by 
\cite[Proposition 6.21]{YZ2} we get
$\opn{Ext}^{i}_{A}(M, R) \cong \opn{Ext}^{i}_{A^{\mrm{op}}}(M, R)$
as bimodules. Hence this bimodule is finite on both sides. We 
conclude that $R$ is weakly bifinite.
\end{proof}

\begin{thm} \label{thm6.1}
Let $A$ be an affine noetherian PI algebra admitting a noetherian 
connected filtration.
\begin{enumerate}
\item
$A$ has a residue complex $\mcal{K}_{A}$.
\item
Let  $B = A / \mfrak{a}$ be a quotient algebra. Then $B$ has a residue
complex $\mcal{K}_{B}$, there is a rigid trace
$\opn{Tr}_{B / A} : \mcal{K}_{B} \to \mcal{K}_{A}$
that is an actual homomorphism of complexes of bimodules, and 
$\opn{Tr}_{B / A} $ induces an isomorphism
\[ \mcal{K}_{B} \cong \opn{Hom}_{A}(B, \mcal{K}_{A}) =
\opn{Hom}_{A^{\mrm{op}}}(B, \mcal{K}_{A}) \subset \mcal{K}_{A} . \]
\end{enumerate}
\end{thm}

\begin{proof}
(1) Is immediate from Proposition \ref{prop6.3} and Theorem 
\ref{thm4.2}. (2) Follows from (1), Theorem \ref{thm5.1} and Corollary 
\ref{cor5.1}.
\end{proof}

Given a set $Z$ of ideals of $A$ we defined a localizing 
subcategory $\msf{M}_{Z} \subset \cat{Mod} A$ in Example 
\ref{exa1.1}. Now let us write $Z^{\mrm{op}}$ for the same set, 
but considered as a set of ideals in the ring $A^{\mrm{op}}$, and let
$\msf{M}_{Z^{\mrm{op}}} \subset \cat{Mod} A^{\mrm{op}}$
be the localizing subcategory. Denote by $\Gamma_{Z}$ and 
$\Gamma_{Z^{\mrm{op}}}$ the two torsion functors respectively.

\begin{cor} \label{cor6.1}
Assume $A$ is like in the theorem, and let $Z$ be a set of 
ideals of $A$. Then
\[ \Gamma_{Z} \mcal{K}_{A} = \Gamma_{Z^{\mrm{op}}} \mcal{K}_{A}
\subset \mcal{K}_{A} . \]
\end{cor}

\begin{proof}
Apply Theorem \ref{thm6.1}(2) to the ideals
$\mfrak{a} = \mfrak{b}_{1} \cdots \mfrak{b}_{n}$ where
$\mfrak{b}_{1}, \ldots, \mfrak{b}_{n} \in Z$,
using formulas (\ref{eqn1.1}) and (\ref{eqn1.2}).
\end{proof}

\begin{cor} \label{cor6.2}
Assume $A$ is like in the theorem.
Let $S$ be a denominator set in $A$, with localization $A_{S}$. Define
$\mcal{K}_{A_{S}} := A_{S} \otimes_{A} \mcal{K}_{A} \otimes_{A} 
A_{S}$.
\begin{enumerate}
\item $\mcal{K}_{A_{S}} \cong A_{S} \otimes_{A} \mcal{K}_{A}
\cong \mcal{K}_{A} \otimes_{A} A_{S}$
as complexes of $A$-bimodules.
\item $\mcal{K}_{A_{S}}$ is a dualizing complex over $A_{S}$.
\end{enumerate}
\end{cor}

\begin{proof}
(1) Let
$Z := \{ \mfrak{p} \in \opn{Spec} A \mid \mfrak{p} \cap S \neq 
\emptyset \}$,
and define Gabriel filters $\mcal{F}_{Z}$ and $\mcal{F}_{S}$ like in 
Examples \ref{exa1.1} and \ref{exa1.3}. According to 
\cite[Theorem VII.3.4]{Ste} these two filters are equal. Hence 
(cf.\ Example \ref{exa1.6}) for each $q$ there is an exact sequence
\[ 0 \to \Gamma_{Z} \mcal{K}^{q}_{A} \to \mcal{K}^{q}_{A} \to
A_{S} \otimes_{A} \mcal{K}^{q}_{A} \to 0 . \]
By symmetry there is an exact sequence
\[ 0 \to \Gamma_{Z^{\mrm{op}}} \mcal{K}^{q}_{A} \to \mcal{K}^{q}_{A} \to
\mcal{K}^{q}_{A} \otimes_{A} A_{S} \to 0 . \]
But by Corollary \ref{cor6.1}, 
$\Gamma_{Z} \mcal{K}^{q}_{A} = 
\Gamma_{Z^{\mrm{op}}} \mcal{K}^{q}_{A}$, which implies
$A_{S} \otimes_{A} \mcal{K}^{q}_{A} \cong 
\mcal{K}^{q}_{A} \otimes_{A} A_{S}$.
Finally use the fact that 
$A_{S} \otimes_{A} A_{S} = A_{S}$.

\medskip \noindent (2) This is proved just like 
\cite[Theorem 1.11(1)]{YZ2} (cf.\ proof of Theorem \ref{thm3.1}(1)).
\end{proof}

Let us remind the reader the definition of a {\em clique} in the prime 
spectrum $\opn{Spec} A$. Suppose $\mfrak{p}$ and $\mfrak{q}$ are 
two prime ideals. If there is a bimodule $M$ that's a subquotient 
of 
$(\mfrak{p} \cap \mfrak{q}) / (\mfrak{p} \mfrak{q})$
and is nonzero torsion-free as $A / \mfrak{p}$-module and as 
$(A / \mfrak{q})^{\mrm{op}}$-module, then we say there is a 
(second layer) link $\mfrak{p} \rightsquigarrow \mfrak{q}$.
The links make $\opn{Spec} A$ into a quiver, and the cliques are 
its connected components. 

\begin{exa} \label{exa6.2}
Suppose $[A : \mbb{K}] < \infty$. The lemma below implies that 
(up to multiplicity of arrows) the link quiver of $A$ coincides 
with the quiver defined by Gabriel in the context of representation 
theory (see \cite{MY}). 
Cliques in this case stand in bijection to blocks of $A$ 
(indecomposable factors), and also to $\opn{Spec} \mrm{Z}(A)$.
\end{exa}

\begin{lem} \label{lem6.1}
Let $A$ be an artinian ring with Jacobson radical $\mfrak{r}$ and 
maximal ideals $\mfrak{p}_{1}, \ldots, \mfrak{p}_{n}$. Then the 
inclusions $\mfrak{r} \subset \mfrak{p}_{i}$ induce an isomorphism 
of $A$-bimodules
\[ \frac{\mfrak{r}}{\mfrak{r}^{2}} \cong 
\bigoplus\nolimits_{i, j} \frac{\mfrak{p}_{i} \cap \mfrak{p}_{j}}
{\mfrak{p}_{i} \mfrak{p}_{j}} . \]
\end{lem}

\begin{proof}
This proof was communicated to us by K. Goodearl. 
Choose orthogonal idempotents $e_{i} \in A$ lifting the central 
idempotents in $A / \mfrak{r}$, so that 
\[ \mfrak{p}_{i} = A (1 - e_{i}) + \mfrak{r} = (1 - e_{i}) A + 
\mfrak{r} . \]
We have
\[ (1 - e_{i})(\mfrak{p}_{i} \cap \mfrak{p}_{j}) = 
(1 - e_{i})\mfrak{p}_{j} \subset 
\mfrak{p}_{i} \mfrak{p}_{j} , \]
and likewise on the right, so each element of 
$(\mfrak{p}_{i} \cap \mfrak{p}_{j}) / (\mfrak{p}_{i} 
\mfrak{p}_{j})$
comes from some element in 
$e_{i} (\mfrak{p}_{i} \cap \mfrak{p}_{j}) e_{j}$. 
But
$e_{i} (\mfrak{p}_{i} \cap \mfrak{p}_{j}) e_{j} =
e_{i} \mfrak{r} e_{j}$. 
We see that the canonical homomorphism
$f: e_{i} \mfrak{r} e_{j} \to 
(\mfrak{p}_{i} \cap \mfrak{p}_{j}) / (\mfrak{p}_{i} \mfrak{p}_{j})$
is surjective. Obviously 
$e_{i} \mfrak{r}^{2} e_{j} \subset \opn{Ker}(f)$. 
On the other hand 
\[ \mfrak{p}_{i} \mfrak{p}_{j} = (1 - e_{i}) A (1 - e_{j}) + 
(1 - e_{i}) \mfrak{r} + \mfrak{r} (1 - e_{j}) + \mfrak{r}^{2} , \]
so 
$e_{i} \mfrak{p}_{i} \mfrak{p}_{j} e_{j} = 
e_{i} \mfrak{r}^{2} e_{j}$. 
Thus 
$\opn{Ker}(f) = e_{i} \mfrak{r}^{2} e_{j}$. 
Finally the isomorphism is obtained by the decomposition
\[ \frac{\mfrak{r}}{\mfrak{r}^{2}} =
\bigoplus_{i,j} e_i \frac{\mfrak{r}}{\mfrak{r}^{2}} e_j =
\bigoplus_{i,j} \frac{e_{i} \mfrak{r} e_{j}}
{e_{i} \mfrak{r}^{2} e_{j}} . \]
\end{proof}

\begin{dfn}
Let $A$ be a noetherian $\mbb{K}$-algebra with an Auslander rigid 
dualizing complex, such that the canonical dimension $\opn{Cdim}$ 
is weakly symmetric. The {\em $q$-skeleton} of $\opn{Spec} A$ is the 
set
\[ \{ \mfrak{p} \in \opn{Spec} A \mid \opn{Cdim} A / \mfrak{p} 
= q \} . \]
\end{dfn}

\begin{prop} \label{prop6.4}
Suppose $A$ is a PI $\mbb{K}$-algebra admitting some noetherian 
connected filtration. Then the $q$-skeleton of $\opn{Spec} A$ is a 
union of cliques.
\end{prop}

\begin{proof}
Within the proof of Theorem \ref{thm4.2} it is shown that 
$\opn{Cdim}$ is constant on cliques in $\opn{Spec} A$.
\end{proof}

\begin{prop} \label{prop6.1}
Suppose $A$ is a prime PI $\mbb{K}$-algebra admitting a noetherian
connected filtration. Let $n := \opn{Cdim} A$ and
$Q := \opn{Frac} A$ the ring of fractions. Then
$\mcal{K}^{-n}_{A} \cong Q$ as $A$-bimodules.
\end{prop}

\begin{proof}
Since $\mcal{K}^{-n}_{A}$ is pure of $\opn{Cdim} = n$ it is 
a torsion-free $A$-module, and  it follows that 
$Q \otimes_{A} \mcal{K}^{-n}_{A} \cong \mcal{K}^{-n}_{A}$.
On the other hand for $q < n$, $\mcal{K}^{-q}_{A}$ is pure of 
$\opn{Cdim} = q < n$, so it is a torsion $A$-module and
$Q \otimes_{A} \mcal{K}^{-q}_{A} = 0$. We see that 
$Q \otimes_{A} \mcal{K}_{A} \cong \mcal{K}^{-n}_{A}[n]$.
But by Theorem \ref{thm6.3} we get
$Q \otimes_{A} \mcal{K}_{A} \cong Q[n]$.
\end{proof}

Given a prime ideal $\mfrak{p}$ in a ring $A$ let 
\[ S(\mfrak{p}) = S_{A}(\mfrak{p}) := \{ a \in A \mid a + \mfrak{p} 
\text{ is regular in } A / \mfrak{p} \} . \]
For a set $Z$ of prime ideals we write
\[ S(Z) := \bigcap_{\mfrak{p} \in Z} S(\mfrak{p}) . \]
According to \cite[Theorem 10 and remarks in Section 3]{Mu},
if $A$ is a noetherian PI affine $\mbb{K}$-algebra and $Z$ is a 
clique of prime ideals in $\opn{Spec} A$ then $S(Z)$
is a denominator set. We get a ring of fractions 
\[ A_{S(Z)} := S(Z)^{-1} \cdot A = A \cdot S(Z)^{-1} . \] 
Furthermore for each 
$\mfrak{p} \in Z$ one has
\[ A / \mfrak{p} \otimes_{A} A_{S(Z)} = Q(\mfrak{p}) = 
\opn{Frac}  A / \mfrak{p} = (A / \mfrak{p})_{S_{A / \mfrak{p}}(0)} 
. \]

For a prime ideal $\mfrak{p}$ denote by $J_{A}(\mfrak{p})$ the
indecomposable injective $A$-module with associated prime
$\mfrak{p}$ (it is unique up to isomorphism). Let 
$\mrm{r}(\mfrak{p})$ denote the Goldie rank of $A / \mfrak{p}$.

We say a clique $Z_{1}$ is a specialization (resp.\ immediate
specialization) of a clique $Z_{0}$ if there exist prime ideals
$\mfrak{p}_{i} \in Z_{i}$ with
$\mfrak{p}_{0} \subset \mfrak{p}_{1}$
(resp.\ and $Z_{i}$ is in the $(q - i)$-skeleton of $\opn{Spec} A$
for some $q$).

Here is the main result of this section.

\begin{thm} \label{thm6.2}
Let $A$ be a PI $\mbb{K}$-algebra admitting a 
noetherian connected filtration, and let $\mcal{K}_{A}$ be its 
residue complex.
\begin{enumerate}
\item For every $q$ there is a canonical $A$-bimodule 
decomposition
\[ \mcal{K}^{-q}_{A} = \bigoplus_{Z} \Gamma_{Z} \mcal{K}^{-q}_{A} \]
where $Z$ runs over the cliques in the $q$-skeleton of
$\opn{Spec} A$.
\item Fix one clique $Z$ in the $q$-skeleton of $\opn{Spec} A$.
Then 
$\Gamma_{Z} \mcal{K}^{-q}_{A} = 
\Gamma_{Z^{\mrm{op}}} \mcal{K}^{-q}_{A}$
is an injective left \tup{(}resp.\ right\tup{)} $A_{S(Z)}$-module, 
and its socle as left \tup{(}resp.\ right\tup{)} $A_{S(Z)}$-module 
is the essential submodule
\[ \bigoplus_{\mfrak{p} \in Z}  \mcal{K}^{-q}_{A / \mfrak{p}} \cong
\bigoplus_{\mfrak{p} \in Z}
\opn{Hom}_{A}(A / \mfrak{p}, \mcal{K}^{-q}_{A})
\subset \Gamma_{Z} \mcal{K}_{A}^{-q} . \]
\item There is a \tup{(}noncanonical\tup{)} decomposition of left 
\tup{(}resp.\ right\tup{)} $A_{S(Z)}$-modules
\[ \Gamma_{Z} \mcal{K}^{-q}_{A} 
\cong \bigoplus_{\mfrak{p} \in Z} 
J_{A}(\mfrak{p})^{\mrm{r}(\mfrak{p})}
\cong \bigoplus_{\mfrak{p} \in Z} 
J_{A^{\mrm{op}}}(\mfrak{p})^{\mrm{r}(\mfrak{p})} . \]
\item $\Gamma_{Z} \mcal{K}^{-q}_{A}$ is an indecomposable
$A$-bimodule.
\item Suppose $Z_{i}$ is a clique in the $(q - i)$-skeleton of
$\opn{Spec} A$, for $i = 0, 1$.
Then $Z_{1}$ is an immediate specialization of $Z_{0}$ iff the
composed homomorphism
\[ \delta_{(Z_{0}, Z_{1})} : \Gamma_{Z_{0}} \mcal{K}^{-q}_{A} \inj
\mcal{K}^{-q}_{A}
\to \mcal{K}^{-q + 1}_{A} \surj \Gamma_{Z_{1}}
\mcal{K}^{-q + 1}_{A} \]
is nonzero.
\end{enumerate}
\end{thm}

\begin{proof}
1.\ Let $I$ be an indecomposable injective module with
associated prime $\mfrak{p} \in Z$, where $Z$ is a clique in the
$q$-skeleton of $\opn{Spec} A$.
Since $A$ satisfies the second layer condition, we get
$\Gamma_{Z} I = I$ and also
$\Gamma_{S(Z)} I = 0$. It follows that
$I \to A_{S(Z)} \otimes_{A} I$
is bijective. Therefore for any other clique $Z'$ in the
$q$-skeleton of $\opn{Spec} A$ we must have
$\Gamma_{Z'} I = 0$.
Because $\mcal{K}^{-q}_{A}$ is an injective module, pure of
dimension $q$, we get the left module decomposition
$\mcal{K}^{-q}_{A} = \bigoplus_{Z} \Gamma_{Z} \mcal{K}^{-q}_{A}$.
By Corollary \ref{cor6.1} and symmetry this is a bimodule
decomposition.

\medskip \noindent
2.\
Clearly $\Gamma_{Z} \mcal{K}^{-q}_{A}$ is an injective
$A_{S(Z)}$-bimodule, and
$\bigoplus_{\mfrak{p} \in Z}  \mcal{K}^{-q}_{A / \mfrak{p}}
\subset \Gamma_{Z} \mcal{K}^{-q}_{A}$
is essential.
Write $Q(\mfrak{p}) := \opn{Fract} A / \mfrak{p}$; then
$\bigoplus_{\mfrak{p}}  \mcal{K}^{-q}_{A / \mfrak{p}}
\cong \bigoplus_{\mfrak{p}}  Q(\mfrak{p})$
is a semi-simple (left and right) $A_{S(Z)}$-module.

\medskip \noindent
3.\ This is because the injective hull of $Q(\mfrak{p})$ as 
$A$-module is $J_{A}(\mfrak{p})^{\mrm{r}(\mfrak{p})}$.

\medskip \noindent
4.\ Assume by contraposition that
$\Gamma_{Z} \mcal{K}^{-q}_{A} = M_{1} \oplus M_{2}$
as bimodules, with $M_{i} \neq 0$. Then the socle
$V = \bigoplus_{\mfrak{p} \in Z}  \mcal{K}^{-q}_{A / \mfrak{p}}$
of $\Gamma_{Z} \mcal{K}^{-q}_{A}$ (as left or right 
$A_{S(Z)}$-module) also decomposes into $V = V_{1} \oplus V_{2}$
with $V_{i} = M_{i} \cap V =
\bigoplus_{\mfrak{p} \in Z_{i}}  \mcal{K}^{-q}_{A / \mfrak{p}}$
and
$Z = Z_{1} \amalg Z_{2}$, $Z_{i} \neq \emptyset$.
Take $\mfrak{p}_{i} \in Z_{i}$ such that there is a second layer
link $\mfrak{p}_{1} \rightsquigarrow \mfrak{p}_{2}$. Recall that
this means there is a bimodule surjection
\[ \mfrak{r} = \frac{ \mfrak{p}_{1} \cap \mfrak{p}_{2} }{
\mfrak{p}_{1} \mfrak{p}_{2} } \surj N \]
with $N$ a nonzero torsion-free module over $A / \mfrak{p}_{1}$ and
$(A / \mfrak{p}_{2})^{\mrm{op}}$.
Then replacing $A$ with $A / \mfrak{p}_{1} \mfrak{p}_{2}$
we retain the link, only now
$Z_{i} = \{ \mfrak{p}_{i} \}$ and
$V_{i} \cong Q(\mfrak{p}_{i})$ as bimodules.

Let $B := A_{S(Z)}$. According to Corollary \ref{cor6.2},
$\mcal{K}_{B} := B \otimes_{A} \mcal{K}_{A} \otimes_{A} B$
is a dualizing complex over $B$. As in the proof of Proposition 
\ref{prop6.1} we get $\mcal{K}_{B}^{-p} = 0$ for all $p < q$, and 
$\mcal{K}_{B}^{-q} = \mcal{K}_{A}^{-q}$.
Hence
$\mcal{K}_{B} = \mcal{K}_{B}^{-q}[q] = \mcal{K}_{A}^{-q}[q]$.

By Lemma \ref{lem3.1} 
\[ \mrm{Z}(B) \cong \opn{End}_{\msf{D}(\cat{Mod} B^{\mrm{e}})}
(\mcal{K}_{B}) \cong
\opn{End}_{B^{\mrm{e}}}(\mcal{K}^{-q}_{B}) \]
as rings. Take 
$\pi \in \opn{End}_{B^{\mrm{e}}}(\mcal{K}^{-q}_{B})$
to be the projection 
$\mcal{K}^{-q}_{B} \surj M_{1}$. So $\pi$ is left multiplication by 
a central idempotent $e \in B$. Since 
$Q(\mfrak{p}_{i}) \cong V_{i} \subset M_{i}$
we see that 
$e \cdot Q(\mfrak{p}_{1}) = Q(\mfrak{p}_{1})$ and 
$e \cdot Q(\mfrak{p}_{2}) = 0$.

Now the bimodule
\[ N_{B} := B \otimes_{A} N \otimes_{A} B \cong
Q(\mfrak{p}_{1}) \otimes_{A} N \otimes_{A} Q(\mfrak{p}_{2}) 
\neq 0 . \]
Being a subquotient of $B$, $e$ centralizes $N_{B}$. We get a 
contradiction
$e \cdot N_{B} = N_{B}$, $N_{B} \cdot e = 0$.

\medskip \noindent
5.\
First assume there is specialization, and choose prime ideals
$\mfrak{p}_{0}, \mfrak{p}_{1}$ as evidence. Then the algebra
$B := A / \mfrak{p}_{0} \otimes_{A} A_{S(Z_{1})}$
is nonzero, having
$Q(\mfrak{p}_{1})$ as a quotient. Thus $B$ is prime.
The complex $\mcal{K}_{B}$, with
$\mcal{K}^{-q}_{B} = \mcal{K}^{-q}_{A / \mfrak{p}_{0}}$,
$\mcal{K}^{-q + 1}_{B} =
\Gamma_{Z_{1}} \mcal{K}^{-q + 1}_{A / \mfrak{p}_{0}}$
and $\mcal{K}^{i}_{B} = 0$ otherwise, is dualizing by Corollary
\ref{cor6.2}. If
$\delta : \mcal{K}^{-q}_{B} \to \mcal{K}^{-q + 1}_{B}$
were zero this would imply that
$B \cong \mrm{H}^{0} \opn{Hom}_{B}(\mcal{K}_{B}, \mcal{K}_{B})$
is decomposable as bimodule, contradicting it being a prime ring.

Conversely assume $\delta_{(Z_{0}, Z_{1})} \neq 0$
and pick some
$\phi \in \Gamma_{Z_{0}} \mcal{K}^{-q}_{A}$
s.t.\
$\delta_{(Z_{0}, Z_{1})}(\phi) \in 
\Gamma_{Z_{1}} \mcal{K}^{-q +1}_{A}$
is nonzero. By part (3) we can find 
$a_{i} \in \mfrak{p}_{1, i} \in Z_{1}$
such that 
\[ 0 \neq \psi = a_{1} \cdots a_{m} \delta_{(Z_{0}, Z_{1})}(\phi) \in 
\mcal{K}^{-q +1}_{A / \mfrak{p}_{1}} \cong Q(\mfrak{p}_{1})  \]
for some $\mfrak{p}_{1} \in Z_{1}$. 
On the other hand there are primes
$\mfrak{p}_{0,1}, \ldots, \mfrak{p}_{0,n} \in Z_{0}$
s.t.\
$\phi \mfrak{p}_{0,1} \cdots \mfrak{p}_{0,n} = 0$,
which implies that
$\psi \mfrak{p}_{0,1} \cdots \mfrak{p}_{0,n} = 0$. 
We conclude that \linebreak
$Q(\mfrak{p}_{1}) \mfrak{p}_{0,1} \cdots \mfrak{p}_{0,m} = 0$,
and therefore
$\mfrak{p}_{0} :=  \mfrak{p}_{0,i} \subset \mfrak{p}_{1}$
for some $i$.
\end{proof}

\begin{exa} \label{exa6.1}
Assume $A$ is finite over $\mbb{K}$,  and let
$A = \prod A_{i}$ be the block decomposition, i.e.\ each $A_{i}$
is an indecomposable bimodule. Then $\opn{Spec} A_{i}$ is a clique
in $\opn{Spec} A$ and 
\[ \mcal{K}_{A} = \mcal{K}_{A}^{0} =
\opn{Hom}_{\mbb{K}}(A, \mbb{K}) = 
\bigoplus_{i} \opn{Hom}_{\mbb{K}}(A_{i}, \mbb{K}) \]
is a decomposition into indecomposable bimodules (cf.\ Example 
\ref{exa6.2}).
\end{exa}

\begin{exa}
Generalizing the previous example, consider a noetherian affine 
$\mbb{K}$-algebra $A$ finite over its center $C$. It is well known 
that $\mfrak{q} \mapsto \mfrak{p} = C \cap \mfrak{q}$
is a bijection from the cliques $Z \subset \opn{Spec} A$ to 
$\opn{Spec} C$ (see \cite[Theorem 11.20]{GW}). 
For $\mfrak{p} \in \opn{Spec} C$ denote by 
$\widehat{C}_{\mfrak{p}}$ the $\mfrak{p}$-adic completion.
The complete semilocal ring
$\widehat{C}_{\mfrak{p}} \otimes_{C} A$
has center $\widehat{C}_{\mfrak{p}}$ and is indecomposable. 
On the other hand, say $\opn{dim} C / \mfrak{p} = n$. Then 
\[ \Gamma_{Z} \mcal{K}_{A}^{-n} = 
\Gamma_{\mfrak{p}} \mcal{K}_{A}^{-n} \cong
\Gamma_{\mfrak{p}} \opn{Hom}_{C}(A, \mcal{K}_{C}^{-n}) \cong
\opn{Hom}_{\widehat{C}_{\mfrak{p}}}(\widehat{C}_{\mfrak{p}} 
\otimes_{C} A, \mcal{K}(\widehat{C}_{\mfrak{p}})) \]
where $\mcal{K}(\widehat{C}_{\mfrak{p}})$ is the dual module from 
Example \ref{exa5.2}. Since 
$\opn{Hom}_{\widehat{C}_{\mfrak{p}}}(-, 
\mcal{K}(\widehat{C}_{\mfrak{p}}))$
is a duality for finite $\widehat{C}_{\mfrak{p}}$-modules we see 
that the indecomposability of the bimodule
$\Gamma_{Z} \mcal{K}_{A}^{-n}$ is equivalent to the 
indecomposability of 
$\widehat{C}_{\mfrak{p}} \otimes_{C} A$.
\end{exa}

\begin{exa} \label{exa6.3}
Consider the PI algebra 
$A = \mbb{K} \bra{x, y} / (y x - q x y, y^{2})$, $q \in \mbb{K}$. 
Assume $\mbb{K}$ is algebraically closed, so the spectrum of $A$ 
consists of the prime ideals $\mfrak{p} := (y)$ and 
$\mfrak{m}_{\lambda} := (y, x - \lambda)$ where 
$\lambda \in \mbb{K}$.
We note that 
\[ \frac{\mfrak{m}_{\lambda} \cap \mfrak{m}_{\mu}}
{\mfrak{m}_{\lambda} \mfrak{m}_{\mu}} \cong
\begin{cases}
\mbb{K} \cdot y \neq 0 & \text{if } \mu = q \lambda \\
0 & \text{otherwise} . 
\end{cases} \]
Thus the cliques are 
$\{ \mfrak{m}_{q^{i} \lambda} \mid i \in \mbb{Z} \}$
and of course $\{ \mfrak{p} \}$.
We see that if $q$ is not a root of unity then we get infinite 
cliques.
\end{exa}

\begin{exa} \label{exa6.4}
Take the quantum plane 
$A := \mbb{K} \bra{x, y} / (y x - q x y)$ with $q$ a primitive 
$l$th root of unity in $\mbb{K}$. The center is 
$C := \mbb{K}[x^{l}, y^{l}]$. Assume $\mbb{K}$ is algebraically 
closed. Let us describe the indecomposable bimodules 
$\Gamma_{Z} \mcal{K}_{A}^{-i}$
and their decomposition into indecomposable left modules
$\Gamma_{Z} \mcal{K}_{A}^{-i} 
\cong \bigoplus_{\mfrak{p} \in Z} 
J_{A}(\mfrak{p})^{\mrm{r}(\mfrak{p})}$.

\medskip \noindent a)
Generically $Q(0)$ is a division ring, and thus
$\mcal{K}_{A}^{-2} \cong J_{A}(0)$. 

\medskip \noindent b)
If $\mfrak{q} \in \opn{Spec} C$ is a curve or a point ($i = 1, 0$)
in the Azumaya locus of $A$ (i.e.\ $x^{l}, y^{l} \notin \mfrak{q}$) 
then $\mfrak{p} := A \mfrak{q}$ is prime. By Tsen's Theorem the 
Brauer group $\opn{Br}(\bsym{k}(\mfrak{q}))$ of the residue field 
$\bsym{k}(\mfrak{q})$ is trivial. Therefore 
$Q(\mfrak{p}) \cong \mrm{M}_{l}(\bsym{k}(\mfrak{q}))$. 
We conclude that $Z := \{ \mfrak{p} \}$ is a clique, 
$\mrm{r}(\mfrak{p}) = l$ and 
$\Gamma_{Z} \mcal{K}_{A}^{-i} \cong
J_{A}(\mfrak{p})^{l}$.

\medskip \noindent c)
If $\mfrak{q} = y^{l} C$  then $\mfrak{p} := y A$ is prime and 
$Q(\mfrak{p})$ is commutative. We conclude that 
$Z := \{ \mfrak{p} \}$ is a clique, 
$\mrm{r}(\mfrak{p}) = 1$ and 
$\Gamma_{Z} \mcal{K}_{A}^{-1} \cong
J_{A}(\mfrak{p})$. Likewise if $\mfrak{q} =  x^{l} C$.

\medskip \noindent d)
If $\mfrak{n} = y^{l} C + (x^{l} - \lambda^{l}) C \in \opn{Spec} C$ 
with $\lambda \neq 0$ then the clique lying above $\mfrak{n}$ is 
$Z := \{ \mfrak{m}_{q^{j} \lambda} \mid j = 0, \ldots, l - 1 \}$;
notation as in the previous example.
The Goldie rank is $\mrm{r}(\mfrak{m}_{q^{j} \lambda}) = 1$ 
and 
$\Gamma_{Z} \mcal{K}_{A}^{0} \cong
\bigoplus_{j = 0}^{l - 1} J_{A}(\mfrak{m}_{q^{j} \lambda})$.
Likewise with $y$ and $x$ interchanged. 

\medskip \noindent e)
Finally the clique lying above 
$\mfrak{n} := y^{l} C + x^{l} C$ is 
$Z := \{ \mfrak{m} \}$ where $\mfrak{m} := xA + yA$, 
$\mrm{r}(\mfrak{m}) = 1$ and 
$\Gamma_{Z} \mcal{K}_{A}^{0} \cong
J_{A}(\mfrak{m})$.
\end{exa}

\begin{que} \label{que6.3}
Are Theorems \ref{thm6.1} and \ref{thm6.2} valid without assuming 
the existence of a noetherian connected filtration? 
\end{que}


\begin{thebibliography}{VdB2}
\bibitem[Aj]{Aj} K.\ Ajitabh, Residue complex for $3$-dimensional
        Sklyanin algebras, Adv.\ in Math.\ \textbf{144} (1999), 
        137-160.
\bibitem[AjSZ]{AjSZ} K.\ Ajitabh, S.P.\ Smith and J.J.\ Zhang,
        Auslander-Gorenstein rings,
        Comm. Algebra \textbf{26} (1998), no. 7, 2159--2180.
\bibitem[ASZ]{ASZ} M. Artin, L.W. Small and J.J. Zhang, 
        Generic flatness for strongly Noetherian algebras,
        J. Algebra \textbf{221} (1999), no.\ 2, 579-610.
\bibitem[ATV]{ATV} M. Artin, J. Tate and M. Van den Bergh, Modules 
        over regular algebras of dimension $3$,
        Invent.\ Math.\ {\bf 106} (1991), 335-388. 
\bibitem[BG]{BG} K.A.\ Brown and K.R.\ Goodearl, Homological aspects of
        Noetherian PI Hopf algebras of irreducible modules and
        maximal dimension,
        J. Algebra \textbf{198} (1997), no. 1, 240--265.
\bibitem[BGK]{BGK} V. Baranovsky, V. Ginzburg and A. Kuznetsov,
        Quiver varieties and noncommutative $\mbf{P}^{2}$,
        preprint. Eprint: math.AG/0103068 at http://arXiv.org.
\bibitem[BM]{BM} G.\ Barou and M.-P.\ Malliavin, Sur la r\'esolution
        injective minimale de l'alg\`ebre enveloppante d'une alg\'ebre
        de Lie r\'esoluble. (French)
        J. Pure Appl. Algebra \textbf{37} (1985), no. 1, 1--25.
\bibitem[Br]{Br} K.A.\ Brown, Fully bounded rings of finite
        injective dimension,
        Quart. J. Math. Oxford Ser., \textbf{41} (1990), 1-13.
\bibitem[Ei]{Ei} D. Eisenbud, ``Commutative Algebra'', Springer, 
        1995.
\bibitem[EG]{EG} P. Etingof and V. Ginzburg, 
        Symplectic reflection algebras, Calogero-Moser space, and 
        deformed Harish-Chandra homomorphism, 
        Invent.\ Math.\  \textbf{147} (2002), 243-348.
\bibitem[GW]{GW} K.R.\ Goodearl and R.B.\ Warfield, Jr.,
        ``An introduction to noncommutative Noetherian rings'',
        London Mathematical Society Student Texts, 16.
        Cambridge University Press, Cambridge-New York, 1989.
\bibitem[KKO]{KKO} A. Kapustin, A. Kuznetsov and D. Orlov,
        Noncommutative Instantons and Twistor Transform,
        Comm.\ Math.\ Phys.\ \textbf{221} (2001), no.\ 2, 
        385-432. 
\bibitem[KS]{KS} M.\ Kashiwara and P.\ Schapira, ``Sheaves
        on Manifolds,'' Springer-Verlag, Berlin, 1990.
\bibitem[Mac]{Mac} S. Mac Lane, Homology. Reprint of the 1975 edition.
        Classics in Mathematics. Springer-Verlag, Berlin, 1995.
\bibitem[Mat]{Mat} H.\ Matsumura, ``Commutative Algebra'', 
        Benjamin/Cummings, Reading, Mass., 1980.
\bibitem[MR]{MR} J.C.\ McConnell and J.C.\ Robson, ``Noncommutative
        Noetherian Rings,'' Wiley, Chichester, 1987.
\bibitem[Mu]{Mu} B.\ M\"{u}ller, Affine noetherian PI-rings have
        enough clans, J.\ Algebra \textbf{97} (1985), 116-129.
\bibitem[MY]{MY} J. Miyachi and A. Yekutieli, Derived Picard
        groups of finite dimensional hereditary algebras, 
        Compositio Math.\ \textbf{129} (2001), no.\ 3, 341-368.
\bibitem[RD]{RD} R.\ Hartshorne, ``Residues and Duality'',
        Lecture Notes in Math.\ \textbf{20},
        Springer-Verlag, Berlin, 1966.
\bibitem[Sp]{Sp} N.\ Spaltenstein, Resolutions of unbounded complexes,
        Compositio Math.\ \textbf{65} (1988), 121-154.
\bibitem[Ste]{Ste} B.\ Stenstr\"{o}m, ``Rings of Quotients,''
        Springer-Verlag, Berlin, 1975.
\bibitem[SV]{SV} J.T. Stafford and M.\ Van den Bergh, 
        Noncommutative curves and noncommutative surfaces, 
        Bull.\ Amer.\ Math.\ Soc.\ \textbf{38} (2001), 171-216. 
\bibitem[SZ]{SZ} J.T.\ Stafford and J.J.\ Zhang, 
        Algebras without noetherian filtrations, preprint.
\bibitem[VdB1]{VdB1} M.\ Van den Bergh, Existence theorems for
        dualizing complexes over non-commutative graded and
        filtered ring,
        J.\ Algebra \textbf{195} (1997), no. 2, 662--679.
\bibitem[VdB2]{VdB2} M.\ Van den Bergh,
        Blowing up of non-commutative smooth surfaces, 
        Mem.\ Amer.\ Math.\ Soc.\ \textbf{154} (2001), 
        no.\ 734, x+140 pp. 
\bibitem[Ye1]{Ye1} A.\ Yekutieli,
        Dualizing complexes over noncommutative graded algebras,
        J.\ Algebra \textbf{153} (1992), 41-84.   
\bibitem[Ye2]{Ye2} A.\ Yekutieli, The residue complex of a
        noncommutative graded algebra, J.\ Algebra \textbf{186} 
        (1996), 522-543.
\bibitem[Ye3]{Ye3} A.\ Yekutieli, 
        Residues and differential operators on schemes,
        Duke Math.\ J.\ \textbf{95} (1998), 305-341.   
\bibitem[Ye4]{Ye4} A.\ Yekutieli, Dualizing complexes, Morita
        equivalence and the derived Picard group of a ring, J. 
        London Math.\ Soc.\ \textbf{60} (1999), 723-746.
\bibitem[Ye5]{Ye5} A.\ Yekutieli, The rigid dualizing complex of a
        universal enveloping algebra,         
        J. Pure Appl.\ Algebra \textbf{150} (2000), 85-93. 
\bibitem[YZ1]{YZ1} A.\ Yekutieli and J.J.\ Zhang, Serre duality for
        noncommutative projective schemes, Proc.\ Amer.\ Math.\ Soc.\
        \textbf{125} (1997), 697-707.
\bibitem[YZ2]{YZ2} A.\ Yekutieli and J.J.\ Zhang,
        Rings with Auslander dualizing complexes,
        J.\ Algebra \textbf{213} (1999), no. 1, 1-51.
\bibitem[YZ3]{YZ3} A.\ Yekutieli and J.J.\ Zhang,
        Dualizing complexes and perverse sheaves on 
        noncommutative ringed schemes, in preparation.
\end{thebibliography}
\end{document}